\documentclass[12pt]{amsart} 
\usepackage{amsmath,amsfonts,amssymb,amsthm,amsxtra,array,amscd}
\usepackage{enumerate}
\usepackage[all]{xy}

\DeclareMathOperator{\om}{om}
\DeclareMathOperator{\ord}{ord}
\DeclareMathOperator{\xt}{xt}
\DeclareMathOperator{\Supp}{Supp}
\DeclareMathOperator{\orb}{orb}

\DeclareMathOperator{\codim}{codim}
\DeclareMathOperator{\dv}{div}
\DeclareMathOperator{\im}{Im}
\DeclareMathOperator{\id}{id}
\DeclareMathOperator{\Coh}{Coh}

\DeclareMathOperator{\Ass}{Ass}
\DeclareMathOperator{\supp}{supp}
\DeclareMathOperator{\Spec}{Spec}
\DeclareMathOperator{\Hom}{Hom}

\DeclareMathOperator{\Coker}{Coker}

\DeclareMathOperator{\reg}{reg}

\newcommand{\lra}{\longrightarrow}

\newcommand{\ka}{{\mathcal A}}
\newcommand{\kc}{{\mathcal C}}
\newcommand{\ko}{{\mathcal O}}
\newcommand{\ke}{{\mathcal E}}
\newcommand{\kf}{{\mathcal F}}
\newcommand{\kg}{{\mathcal G}}
\newcommand{\kh}{{\mathcal H}}
\newcommand{\ki}{{\mathcal I}}
\newcommand{\kl}{{\mathcal L}}
\newcommand{\km}{{\mathcal M}}
\newcommand{\kn}{{\mathcal N}}
\newcommand{\kq}{{\mathcal Q}}
\newcommand{\kt}{{\mathcal T}}

\newcommand{\fm}{\mathfrak{m}}
\newcommand{\fp}{\mathfrak{p}}

\newcommand{\rays}{{\Delta(1)}}
\newcommand{\N}{{\mathbb N}}
\newcommand{\A}{{\mathbb A}}
\newcommand{\Z}{{\mathbb Z}}
\newcommand{\R}{{\mathbb R}}

\renewcommand{\P}{{\mathbb P}}

\newcommand{\E}{{\mathbb E}}

\newcommand{\quot}{{//}}

\newcommand{\comment}[1]{}

\begin{document}

\newtheorem{sub}{}[section]
\newtheorem{subsub}{}[sub]
\newtheorem{subsubsub}{}[sub]

\title[Equivariant Primary Decomposition]{Equivariant Primary 
Decomposition and Toric Sheaves}

\author[Perling]{M.~Perling} 
\address{Ruhr-Universit\"at Bochum\\
Fakult\"at f\"ur Mathematik
\newline Universit\"atsstrasse 150\\
\newline D-44780 Bochum}
\email{Markus.Perling@rub.de}
  
\author[Trautmann]{G.~Trautmann}
\address{Universit\"at Kaiserslautern\\
  Fachbereich Mathematik\\
\newline  Erwin-Schr\"odinger-Stra{\ss}e
\newline D-67663 Kaiserslautern}
\email{trm@mathematik.uni-kl.de}

\footnotetext[1]{Research of both authors was supported by the DFG
Schwerpunktprogramm 1094}


\begin{abstract}
We study global primary decompositions in the category of sheaves on a scheme
which are equivariant under the action of an algebraic group. We show that
equivariant primary decompositions exist if the group is connected.
As main application we consider the case of varieties which are quotients
of a quasi-affine variety by the action of a diagonalizable group and
thus admit a homogeneous coordinate ring, such as toric varieties.
Comparing these decompositions with primary decompositions of 
graded modules over the homogeneous coordinate ring,
we show that these are equivalent if the action of the diagonalizable
group is free. We give some specific examples for the case of toric varieties.

\end{abstract}
\vspace{2cm}
\maketitle

\vspace{1cm}

\tableofcontents

\newpage

\section*{Introduction}

Global primary decompositions of coherent sheaves can be found
in \cite[\S 3]{EGAIV} in the case of schemes, and in the analytic case their 
existence was proved by Y. T. Siu in \cite{Siu}, see also \cite{Cy}.
In this paper we study global primary decomposition of equivariant
coherent sheaves on algebraic varieties and schemes with a group action,
and in particular of coherent toric sheaves on toric varieties.
Equivariant primary decomposition is equivalent to primary decomposition
of coherent sheaves on the quotient stack $[X/G]$ of a scheme $X$ by
an algebraic group $G$. From a more geometric point of view it is natural
to ask whether an equivariant primary decomposition gives rise to a primary
decomposition of sheaves on a good quotient $X//G$. This is not true in
general and we will study this fact to some extend for the case of sheaves
on varieties which admit homogeneous coordinates, in particular on toric
varieties.

Following \cite{BerchtoldHausen}, we say that a variety $X$ admits homogeneous
coordinates if there exists an affine
variety $Z$ together with an action of a diagonalizable algebraic
group $H$ and an open subset $W$ such that $X$ is a good quotient
of $W$ by $H$. Then the coordinate ring $S$ of $Z$ acquires a grading
by the character group of $H$ and $S$ serves
as a homogeneous coordinate ring for $X$ with respect to this grading.
This setting comes with a natural sheafification functor $F \mapsto 
\widetilde{F}$,
which maps a graded $S$-modules to a quasi-coherent sheaf on $X$.
This generalizes the usual homogeneous coordinate
rings for projective spaces and toric varities.

Naively, for primary decomposition of sheaves on $X$ it would seem to be 
sufficient
to look at primary decompositions of graded $S$-modules.
However, a priori it is not clear (and actually not true, see examples
\ref{singexample1} and \ref{singexample2}) that a graded primary
decomposition of some $S$-module $F$ yields
a proper primary decomposition of sheaves of $\widetilde{F}$. We show in 
Theorem
\ref{tordec2} that this at least holds if $X$ is a geometric quotient of $W$ 
by $H$.

In section \ref{supplements} we review some general facts
on equivariant quasi-coherent sheaves, together with some clarifying remarks. 
Most of this material should be well-known.
In section \ref{locco} we transport material from \cite{SiuTrm} on gap
sheaves and local cohomology to the equivariant setting.
In section \ref{sectepd} we give an essentially complete account on the theory
of equivariant primary decomposition on general locally noetherian schemes
over a base field of characteristic zero; this is exemplified in section
\ref{ptwo} for some examples of torus-equivariant sheaves on $\P_2$.
In section \ref{sqp} we consider varieties which admit homogeneous
coordinates and study to some extend the descent of $G$-equivariant
sheaves on $Z$. Here, $G$ is an algebraic group which contains
$H$ as a normal subgroup.
We obtain our main results in section
\ref{quotientprimdec}, where we compare graded $G$-equivariant
primary decomposition over $S$ with $G/H$-equivariant primary
decomposition over $X$. Finally, in section \ref{toricsh}, as an
explicit application, we will construct equivariant primary decompositions for
sheaves of Zariski differentials over toric varieties.

{\bf Acknowledgements.} We want to thank the referee for some interesting
remarks which enabled us to generalize our results to a much wider context.

\section{Supplements on equivariant sheaves}\label{supplements}
All schemes in this paper are supposed to be locally of finite type over some
fixed algebraically closed field $k$ of characteristic $0$. Let $G$
be a linear algebraic group over $k$ and let
\[
G \times X \overset{\sigma}{\lra} X
\]
be an action on the scheme $X$.   A coherent $\ko_X$--module $\ke$ is called
$G$--linearized or $G$--equivariant, if there is an isomorphism
\[
\sigma^\ast \ke \overset{\Phi}{\lra} p_2^\ast \ke
\]
such that the following cocycle condition is satisfied. 

Letting $\mu$ denote the multiplication morphism of $G$ and $p_2$, resp.\
$p_{23}$, the projections 
$G \times X \to X$, resp.\ $G \times G \times X \to
G \times X$, then the diagram (modulo canonical isomorphisms)
\[
\UseComputerModernTips
\xymatrix{
\bigl(\sigma \circ (\mu \times 1)\bigr)^\ast \ke \ar[r]^\approx \ar[d]_{(\mu
  \times 1)^\ast \Phi} & \bigl(\sigma
\circ (1 \times \sigma)\bigr)^\ast\ke \ar[r]^{(1 \times \sigma)^\ast \Phi} &
  \bigl(p_2 \circ(1 \times \sigma )\bigr)^\ast \ke \ar[d]^\approx \\
\bigl(p_2 \circ (\mu \times 1)\bigr)^\ast \ke \ar[r]^\approx & (p_2 \circ
  p_{23})^\ast \ke & (\sigma \circ p_{23})^\ast \ke \ar[l]^{p_{23}^\ast \Phi}
}
\]
is supposed to be a commutative diagram of isomorphisms.   Here we follow the
notation of \cite{MumFog}. An isomorphism of this type will be called
a  ($G$--)linearization or a ($G$--)equivariance. More precisely,
an equivariant sheaf is a pair $(\ke, \Phi)$ of
a sheaf $\ke$ and an ($G$--)equivariant structure $\Phi$ on $\ke$. By abuse of 
notation we will simply call $\ke$ an ($G$--)equivariant sheaf, when there is 
no doubt about  $\Phi$.

When $g\in G(k)$ is a closed point, we  are given
an isomorphism $X \xrightarrow{g} X$ and the immersion 
$X\overset{i_g}{\hookrightarrow} G\times X$ defined by $x \mapsto (g,x)$. 
The restriction of an equivariance  $\Phi$ of $\ke$  to $X$ 
induces an isomorphism $\Phi_g$ via
\[
\UseComputerModernTips
\xymatrix{
i^\ast_g \sigma^\ast \ke \ar[r]^{i^\ast_g \Phi} \ar[d]_\approx & i_g^\ast
p_2^\ast \ke\ar[d]_\approx \\
g^\ast \ke \ar[r]^{\Phi_g} & \ke\;.
}
\]
The family of these isomorphisms then satisfies the following cocycle
condition, i.e. the commutativity of the diagram
$$
\UseComputerModernTips
\xymatrix{
g^\ast_2g^\ast_1 \ke \ar[d]_\approx \ar[r]^{g_2^\ast \Phi_{g_1}} & g_2^\ast\ke
\ar[d]_\approx^{\Phi_{g_2}} \\
(g_1g_2)^\ast\ke \ar[r]_{\Phi_{g_1g_2}} & \ke
}\eqno(A)
$$
for any two elements $g_1, g_2 \in G(k)$.   
To see this, one can restrict the
diagram on $G \times G \times X$ to a fibre $\bigl\{(g_1, g_2)\bigr\}
\times X$ to obtain a diagram (A).   Conversely, this conditon implies
the original condition on $\Phi$, by considering the
homomorphisms 
\begin{gather*}
\alpha = (p_{23}^\ast \Phi) \circ (1 \times \sigma)^\ast \Phi\quad \text{ and }
\quad\beta = (\mu \times 1)^\ast \Phi\,,\\
\begin{array}{ccc}
\bigl(\sigma \circ (1 \times \sigma)\bigr)^\ast \ke & & (p_2 \circ
p_{23})^\ast \ke\,.\\
& \raisebox{3.25ex}[-3.25ex]{$\overset{\alpha}{\lra}$} & \\
& \raisebox{4.75ex}[-4.75ex]{$\underset{\beta}{\lra}$} & \\
\end{array}
\end{gather*}
Then commutativity of the diagrams (A) implies that their restrictions to
each fibre $\bigl\{(g_1, g_2) \bigr\} \times X$ coincide.   Then $\alpha =
\beta $ follows from Lemma \ref{restr} below.
\vskip5mm

\begin{sub}{\bf Remark:}\label{existeq}\; \rm
An equivariance $\Phi$ is uniquely determined by the cocycle $(\Phi_g)$.
One should note however that a given  system  $(\varphi_g)$ of isomorphisms
$g^\ast\ke\to\ke$ with (A), called cocycle, need
not be derived from an isomorphism $\Phi$. To see that, let first $(\Phi_g)$
be derived, and let $G\overset{\chi}{\lra}k^\ast $ be a group homomorphism
which is not a character. Such homomorphisms exist for algebraically closed 
fields. Then the family $(\chi(g)\Phi_g)$ is again a 
cocycle. If that should be derived from an equivariance $\Psi$ with
$\Psi_g=\chi(g)\Phi_g$ for any $g$, then by uniqueness $\Psi=\chi\Phi$. Then
$\chi$ would be a homothetic isomorphism of $\sigma^\ast\ke$ or 
$p_2^\ast\ke$, and then a character. On the other hand, when a cocycle
$(\varphi_g)$  is given, one could ask for a group homomorphism
$G\overset{\chi}{\lra}k^\ast $ such that $(\chi(g)\varphi_g)$ is derived
from an equivariance.
\end{sub}

In \cite{kly} A.~Klyachko showed that on a complete
toric variety a vector bundle is already equivariant if only $t^\ast E\simeq E$
for any torus transformation $t$. 
The following proposition shows an analogous statement, namely that under the 
same
conditions an equivariance exists
in case of invariant subsheaves or quotient sheaves of an equivariant sheaf.

Any action $\sigma$ defines an isomorphism
$G \times X \overset{\tilde{\sigma}}{\lra} G \times X$
as the morphism $(p_1, \sigma)$ with underlying map $(g,x) \to (g,gx)$ for 
closed points,
such that $\sigma = p_2 \circ \tilde{\sigma}$.   It follows that $\sigma$ is a
\textbf{flat} morphism, as is $p_2$, see \cite[1.4]{thom}.   Therefore, for any
submodule $\kf \subset \ke$, there are the pulled back submodules $\sigma^\ast
\kf \subset \sigma^\ast \ke$ and $p_2^\ast \kf \subset p_2^\ast \ke$.
\vskip5mm

\begin{sub}{\bf Proposition:}\label{inv1}
  Let $\ke$ be a $G$--equivariant coherent sheaf on $X$ under the action $G
  \times X \to X$ with isomorphism $\Phi$, and let $\kf \subset \ke$ be a
  coherent submodule.   Suppose that for any $g \in G(k)$ the isomorphism
  $\Phi_g$ induces an isomorphism $g^\ast \kf \to \kf$ of the subsheaves.
    Then $\Phi$ induces an isomorphism $\Psi$ for $\kf$, by which $\kf$
  becomes an equivariant (sub)sheaf as indicated by the diagram
\[
\begin{array}{ccc}
\sigma^\ast \ke & \overset{\Phi}{\underset{\approx}{\lra}} & p_2^\ast \ke\\
\cup & & \cup\\
\sigma^\ast \kf & \overset{\Psi}{\underset{\approx}{\lra}} & p_2^\ast \kf\,.
\end{array}
\]
The analogous statement is true for coherent quotients of $\ke$.
\end{sub}
\begin{proof}
Consider the diagram
\[
\UseComputerModernTips
\xymatrix{
0 \ar[r] &  \sigma^\ast \kf \ar[drr]^>>>>>\alpha\ar[r] & \sigma^\ast\ke
\ar[r]\ar[d]^\Phi & \sigma^\ast (\ke/\kf) \ar[r] & 0\\
0 \ar[r] & p_2^\ast\kf \ar[r] &  p_2^\ast\ke \ar[r] &  p_2^\ast (\ke/\kf) 
\ar[r]& 0 
}
\]
and let $\alpha$ denote the composition.   We are going to show that $\alpha =
0$.   Then $\Phi$ induces a homomorphism $\Psi$ as in the statement, which
must be an isomorphism.   By assumption, for any $g \in G(k)$ the restricted
diagram becomes
\[
\UseComputerModernTips
\xymatrix{
0 \ar[r] & g^\ast\kf \ar@{-->}[d]^\approx \ar[r] & g^\ast \ke
\ar[d]_\approx^{\Phi_g} \ar[r] & g^\ast (\ke/\kf) \ar[r] & 0\\
0 \ar[r] & \kf \ar[r] & \ke \ar[r] & \ke/\kf \ar[r] & 0
}
\]
with induced isomorphism for $\kf$.   It follows that the restriction
$\alpha_g = 0$ for any $g$.   By the following Lemma \ref{restr} and the
fact that in characteristic zero $G$ always is reduced, $\alpha =
0$.   The induced homomorphism $\sigma^\ast \kf \xrightarrow{\Psi} p_2^\ast
\ke$ is automatically injective.   If $\kc$ denotes its cokernel, then our
assumption implies also that $\kc_g = i_g^\ast \kc = 0$ for any $g$, hence
$\kc = 0$.

Finally, the cocycle condition for $\ke$ implies that for $\kf$, because all
the involved pull backs are submodules of the corresponding pull backs of
$\ke$.
\end{proof}
\vskip3mm

\begin{sub}{\bf Corollary:}\label{invsub}
  Let $G \times X \xrightarrow{\sigma} X$ be an action as above, and let $Y
  \subset X$ be a subscheme, which is invariant under each $g \in G(k)$, i.e.\
  the isomorphism $g$ of $X$ induces an isomorphism of $Y$.   Then $\ko_Y$ is
  equivariant or $\sigma$ induces an action $G \times Y \to Y$ on $Y$.
\end{sub}
\begin{proof}
  It is canonical to verify that $G \times Y \xrightarrow{\sigma} X$
  factorizes through $Y$ if and only if $\ko_Y$ is a $G$--equivariant quotient
  of $\ko_X$.   The assumption implies that $\sigma^\ast \ko_X
  \xrightarrow{\Phi} p_2^\ast \ko_X$ induces isomorphisms $g^\ast \ko_Y \to
  \ko_Y$ for all $g \in G(k)$.   By Proposition \ref{inv1}, $\Phi$ induces an
  isomorphism $\sigma^\ast \ko_Y \xrightarrow{\Psi} p_2^\ast \ko_Y$, by which
  $\ko_Y$ is $G$--equivariant.
\end{proof}
\vskip3mm

\begin{sub}{\bf Lemma:}\label{restr}
  Let $S$ be a reduced scheme, $\km$ a coherent module on $S \times X$, $\kn$
  a coherent module on $X$ and let $\km \xrightarrow{\alpha} p_2^\ast \kn$ be
  a homomorphism.   If the composition 
\[
\alpha_s : \km \overset{\alpha}{\lra}   p_2^\ast \kn \lra i_s^\ast p_2^\ast \kn
\cong \kn
\]
of $\alpha$ with the restriction map to each fibre $\{s\} \times X$, $s \in
S(k)$, is zero, then $\alpha = 0$.
\end{sub}

\begin{proof}
  It is sufficient to consider an affine neighbourhood of an arbitrary closed
  point of $S \times X$.   Therefore, we may assume that both $S$ and $X$ are
  affine.   It is then sufficient to prove that $\Gamma(\alpha) = 0$.   
We are now given the diagrams
\[
\UseComputerModernTips
\xymatrix{\Gamma(S \times X,\, \km)\ar[r]^{\Gamma(\alpha)}
\ar[dr]_{\Gamma(\alpha_s)} &
\Gamma (S \times X,  p_2^\ast \kn)\ar[d]\ar[r]^{\approx} &
\Gamma(S,\ko_S) \otimes_k \Gamma(X,\kn)\ar[d] & \\
 & \Gamma(X_s,\; i_s^\ast p_2^\ast\kn)\ar[r]^{\approx} & 
k(s)\otimes_k\Gamma(X,\kn)\ar[r]^{\approx} & \Gamma(X,\kn)
}
\]
Let $f \in \im \Gamma(\alpha)$.   We can write
\[
f = \sum^n_{\nu=1} a_\nu \otimes f_\nu
\]
with $a_\nu \in \Gamma(S,\ko_S)$ and $f_\nu \in \Gamma(X,\kn)$, and we may
assume that $a_1, \dots, a_n$ are linearly independent over $k$.
By assumption, the evaluation at any $s \in S(k)$ yields
\[
\sum_\nu a_\nu(s) f_\nu = 0\,.
\]
Because of the linear independence, we can find points $s_1, \dots, s_n \in
S(k)$, such that the matrix $\bigl(a_\nu(s_\mu)\bigr)$ is invertible.   It
follows that $f_1 = \cdots = f_n = 0$.
\end{proof}
\vskip3mm

\begin{sub}\rm {\bf Equivariant homomorphisms}\label{equihom}\\
Let $(\ke, \Phi)$ and $(\kf, \Psi)$ be $G$--equivariant coherent sheaves on $X$
with respect to an action $G \times X \xrightarrow{\sigma} X$.   A
homomorphism $\ke \xrightarrow{a} \kf$ is called $G$--invariant if it is
compatible with the isomorphisms $\Phi$ and $\Psi$, i.e.\ if the diagram
\[
\UseComputerModernTips
\xymatrix{
\sigma^\ast \ke \ar[r]^{\sigma^\ast a} \ar[d]^\Phi_\approx & \sigma^\ast\kf 
\ar[d]^\Psi_\approx\\
p^\ast_2\ke \ar[r]^{p_2^\ast a} & p_2^\ast\kf
}
\]
is commutative.  As in the case of subsheaves, it is sufficient that
$a$ is compatible with $\Phi_g$ and $\Psi_g$ for any $g \in G(k)$, i.e.\ if
the diagrams
\[
\UseComputerModernTips
\xymatrix{
g^\ast\ke \ar[r]^{g^\ast a} \ar[d]^{\Phi_g}_\approx & g^\ast \kf
\ar[d]^{\Psi_g}_\approx\\
\ke  \ar[r]^a & \kf
}
\]
are commutative for any $g \in G(k)$.  In order to prove this, let $\ka$ be
the image of $a$.   Then we have the commutative diagrams
\[
\UseComputerModernTips
\xymatrix{
g^\ast \ke \ar@{->>}[r]^{g^\ast a} \ar[d]^{\Phi_g}_\approx & g^\ast\ka\;\;
\ar@{>->}[r] \ar@{-->}[d]^{\alpha_g}_\approx & g^\ast \kf \ar[d]^{\Psi_g}\\
\ke \ar@{->>}[r] & \ka\;\; \ar@{>->}[r] & \kf
}
\]
and induced isomorphisms $\alpha_g$.   Then Proposition \ref{inv1} yields the 
claim.   
\end{sub}
\vskip3mm

\begin{sub}\rm {\bf Remark:}\label{abcat}\;
The category $\Coh^G(X)$ of $G$--equivariant coherent sheaves with 
$G$--invariant
homomorphisms is an abelian category, see also \cite[1.2, 1.4]{thom}.
\end{sub}
\vskip3mm

\begin{sub}\rm {\bf Dual actions:}\label{dact}\;
Let $G$ be a linear algebraic group over $k$ and let $G\times
X\xrightarrow{\sigma} X$ be an action on an affine variety over $k$. Denoting
by $S=k[X]$ the coordinate ring, $\sigma$ and the multiplication map $\mu$ of
$G$ induce dual actions
\[
S\xrightarrow{\sigma^\ast} k[G]\underset{k}{\otimes} S\ \text{ and }\
k[G]\xrightarrow{\mu^\ast} k[G]\underset{k}{\otimes} k[G]\ .
\]
A dual action on an $S$--module $M$ is a $k$-linear map
\[
M\xrightarrow{\vartheta} k[G]\underset{k}{\otimes} M
\]
such that $(\mu^\ast\otimes 1)\circ \vartheta=(1\otimes \vartheta)\circ
\vartheta$ and $(\varepsilon\otimes 1)\circ\vartheta=\id_M$ where
$k[G]\xrightarrow{\varepsilon} k$ is the evaluation at $e\in G$ and such that 
\[
\vartheta(sm)=\sigma^\ast(s)\vartheta (m)
\]
for all $s\in S$ and $m\in M$, see \cite{MumFog}. If $\varepsilon_g$ denotes the
evaluation map at $g\in G(k)$, then a dual action induces homomorphisms
\[
M\xrightarrow{\vartheta_g} M\ \text{ and } S\xrightarrow{\sigma^\ast_g} S
\]
with $\vartheta_g=(\varepsilon _g\otimes 1)\circ \vartheta$ and
$\sigma^\ast_g=(\varepsilon_g\otimes 1)\circ\sigma^\ast$, satisfying the rule
\[
\vartheta_g(sm)=\sigma^\ast_g(s)\vartheta_g(m)\ .
\]
The compatibility of $\mu^\ast$ and $\vartheta$ implies that
\[
\vartheta_{g_1g_2}=\vartheta_{g_2}\circ\vartheta_{g_1}\ \text{ and }
\sigma^\ast_{g_1g_2}=\sigma^\ast_{g_2}\circ \sigma^\ast_{g_1}
\]
for any two elements of $G(k)$. In particular, $\vartheta$ defines an action
of $G$ on $M$, and each $\vartheta_g$ is an isomorphism. Let now
\[
(k[G]\underset{k}{\otimes} S)\otimes^\sigma_S M\xrightarrow{\varphi}
k[G]\otimes_k M
\]
be defined by $\varphi(a\otimes s\otimes m)=(a\otimes s)\ .\ \vartheta(m)$
where $\otimes^\sigma_S$ indicates that the tensorproduct is taken with
respect to $\sigma^\ast$, and where the right hand side is the scalar
multiplication of the $k[G]\underset{k}{\otimes} S$ -- module $k[G]\otimes_k
M$.

Then $\varphi$ is $k[G]\otimes_kS$--linear and defines a sheaf homomorphism
\[
\sigma^\ast \widetilde{M}\xrightarrow{\widetilde{\varphi}}
p^\ast_2\widetilde{M}
\]
of the pull backs of the associated sheaf $\widetilde{M}$.

In order to see that $\widetilde{\varphi}$ is a $G$--linearization we can as
well define $S$--homomorphisms $S\otimes^{\sigma_g}_S M\xrightarrow{\varphi_g}
M$ by $\varphi_g(s\otimes m)=s\vartheta_g(m)$ for $g\in G(k)$.

Note here that $\varphi_g(1\otimes
sm)=\vartheta_g(sm)=\sigma^\ast_g(s)\vartheta_g(m)=\vartheta_g(\sigma^\ast_g(s)
\otimes
m)$ and that $\varphi_g$ is an $S$--isomorphism. The maps $\varphi_g$ are
obtained from $\varphi$ by evaluation at $g$ and thus the induced sheaf
homomorphisms
\[
g^\ast \widetilde{M}\xrightarrow{\widetilde{\varphi}_g} \widetilde{M}
\]
are restrictions of $\widetilde{\varphi}$ to $\{g\}\times X$. It follows that
$\widetilde{\varphi}$ is an isomorphism. Moreover, the rule
$\vartheta_{g_1g_2}=\vartheta_{g_2}\circ\vartheta_{g_1}$ implies that the
cocycle condition holds for $\widetilde{\varphi}$. Thus, a dual action
$\vartheta$ on $M$ gives rise to a $G$--linearization
$\Phi_\vartheta:=\widetilde{\varphi}$ on $\widetilde{M}$. Since
conversely any $G$--linearization $\Phi$ on $\widetilde{M}$ defines a dual
action on $M$ as composition 
\[
\Gamma(X, \widetilde{M})\to\Gamma(G\times X, \sigma^\ast
\widetilde{M})\xrightarrow{\Gamma(\Phi)}\Gamma(G\times X, p_2^\ast
\widetilde{M})\cong k[G]\otimes \Gamma(X, \widetilde{M})
\]
and since this is inverse to $\vartheta\mapsto\Phi_\vartheta$, we have the
\end{sub}

\begin{sub}{\bf Proposition:}\label{dactpr} The map $\vartheta\to\Phi_\vartheta$ 
defines a
$1-1$ correspondence between dual actions of the $S$--module $M$ and the
$G$--linearizations of $\widetilde{M}$.
\end{sub}
\vskip4mm

\begin{sub}{\bf The case of a diagonalizable group:}\label{dactdiag}\rm\quad  
Let  
$G\times X\xrightarrow{\sigma} X$ and $M$ be as in the previous section
and assume that $G$ is diagonalizable with the character group $X(G)$ being
a basis of the $k$--algebra $k[G]$. Given a dual action
$M\xrightarrow{\vartheta} k[G]\otimes M$ of the $S$--module $M$, we obtain an
$X(G)$--grading of $M$ by the subgroups
\[
\begin{array}{lcl}
M_\chi & = & \{m\in M\ |\ \vartheta(m)=\chi\otimes m\}\\
       & = & \{m\in M\ |\ \vartheta_g(m)=\chi(g)m\ \text{ for all } g\in G\}
\end{array}
\]
for any character $\chi$, and such that $S_{\chi_1} M_{\chi_2}\subset
M_{\chi_1+\chi_2}$, where also $S$ is $X(G)$-graded via $\sigma^\ast$.
Conversely, given an $X(G)$--grading of $M$, one can
define a dual action $\vartheta$ via $\vartheta(m_\chi):=\chi\otimes m_\chi$
for $\chi$--homogeneous elements. We obtain the
\end{sub}

\begin{sub}{\bf Corollary:}\label{dactdiagcor} Let $G$ be a diagonalizable 
group
  with an action
  $G\times X\xrightarrow{\sigma} X$ on an affine variety, and let $M$ be a
  $k[X]$--module. There is a $1:1$ correspondence between $X(G)$--gradings of
  $M$ and $G$--linearizations of the associated sheaf $\widetilde{M}$ on $X$.
\end{sub}
\vskip4mm

\begin{sub}\rm {\bf Equivariance of divisors and divisorial sheaves:}
\label{equidiv}\\
To any Weil--divisor $D =\sum n_Y Y$ on a normal variety $X$ one can associate
a reflexive sheaf $\ko_X(D)$ of rank 1 as the sheaf $i_\ast
\ko_{X_{\reg}}(D|X_{\reg})$ where $i$ denotes the inclusion of $X_{\reg}$ in
$X$.
Alternatively, for any open $U$ we have
$$\ko_X(D)(U)=\{f\in k(X) | \ord_Y(f)
+ n_Y\geq 0 \text{ for any component } Y \text{ of } D \text{ with }Y\cap
U\neq \emptyset\}.$$

The map $D\mapsto \ko_X(D)$ induces a one-to-one correspondence between the
divisor class group $Cl(X)$ and the set of isomorphism classes of reflexive
rank--1 sheaves on $X$, see e.g. \cite{Re}, and it reduces to the usual
correspondence for the Picard group if $X$ is smooth.

If $Y\xrightarrow{\varphi} X$ is a flat morphism of fixed relative dimension,
then $\varphi^\ast \ko_X(D)$ and $\ko_Y(\varphi^\ast D)$ are canonically
isomorphic.

Let now $G\times X\xrightarrow{\sigma} X$ be an action of a linear algebraic
group on a normal variety over $k$. A Weil--divisor $D$ on $X$ is then called
invariant if all of its components are invariant. If so,
\[
\sigma^\ast D=\sum n_Y(G\times Y)= p_2^\ast D,
\]
and then $\ko_X(D)$ is equipped with a canonical equivariance $\Phi$ defined
by the diagram
\[
\UseComputerModernTips
\xymatrix{
\sigma^\ast\ko_X(D)\ar[r]^\Phi\ar[d]^\approx &
p_2^\ast\ko_X(D)\ar[d]^\approx\\
\ko_{G\times X}(\sigma^\ast D)\ar@{=}[r] & \ko_{G\times X}(p_2^\ast D)
}\raisebox{-8ex}{,}
\]
where the vertical arrows are the canonical isomorphisms. If $D$ is in
addition effective, then the canonical section $\ko_X\to \ko_X(D)$ is
$G$--invariant with respect to $\Phi$.
\end{sub}
\vskip3mm

\begin{sub}\rm {\bf Equivariance of composed sheaves}\label{equiv}\\
Let $Y \xrightarrow{f} X$ be a morphism of schemes/$k$ and let $\kf$ and $\kg$
be coherent $\ko_X$--modules.   There are canonical homomorphisms
\begin{align*}
f^\ast \kf \otimes_{\ko_Y} f^\ast \kg & \underset{\approx}{\lra} f^\ast (\kf
\otimes_{\ko_X} \kg)\\
\intertext{and}
f^\ast \kh\!\om_{\ko_X} (\kf, \kg)& \lra \kh\!\om_{\ko_Y} (f^\ast \kf, f^\ast
\kg)\,. 
\end{align*}
The first is always an isomorphism, the latter is an isomorphism if $f$ is
\textbf{flat}, see \cite[Ch.\ 0, 4.3, 4.4]{EGAI}.

For flat $f$ there are also canonical isomorphisms
\[
f^\ast \ke\!\xt^i_{\ko_X} (\kf,\kg) \lra \ke\!\xt^i_{\ko_Y} (f^\ast \kf,
f^\ast \kg)\,,
\]
see \cite[Ch.\ 0, 12.3.4]{EGAIII}.
\end{sub}
\vskip3mm

\begin{sub}{\bf Lemma:}\label{indequi}
  Let $\kf$ and $\kg$ be $G$--equivariant coherent sheaves on $X$ with respect
  to an action $G \times X \xrightarrow{\sigma} X$.  Then the sheaves $\kf
  \otimes_{\ko_X} \kg$ and $\ke\!\xt^i_{\ko_X} (\kf,\kg)$ are also
  $G$--equivariant. 
\end{sub}
\begin{proof}
  Let $\sigma^\ast \kf \xrightarrow{\Phi} p_2^\ast \kf$ and $\sigma^\ast \kg
  \xrightarrow{\Psi} p_2^\ast \kg$ be the isomorphisms of equivariance for
  $\kf$ and $\kg$.   We obtain the canonically commutative diagrams
\[
\UseComputerModernTips
\xymatrix{
\sigma^\ast \kf \otimes \sigma^\ast \kg \ar[d]^{\Phi \otimes \Psi}_\approx
\ar[r]_\approx & \sigma^\ast(\kf \otimes \kg)\ar@{-->}[d]_\approx\\
p_2^\ast \kf \otimes p_2^\ast \kg \ar[r]_\approx & p_2^\ast(\kf \otimes \kg)
}
\]
and
\[
\UseComputerModernTips
\xymatrix{
\sigma^\ast \ke\!\xt^i(\kf,\kg) \ar@{-->}[d]_\approx \ar[r]_\approx &
\ke\!\xt^i(\sigma^\ast \kf, \sigma^\ast
\kg)\ar[d]_\approx^{\varepsilon^i(\Phi,\Psi)}\\
p_2^\ast \ke\!\xt^i(\kf,\kg) \ar[r]_\approx & \ke\!\xt^i(p^\ast_2\kf, p_2^\ast
\kg)\,,
}
\]
where $\varepsilon^i(\Phi,\Psi)$ is induced by $\Phi$ and $\Psi$, and where
the horizontal homomorphisms are isomorphisms due to flatness of $\sigma$ and
$p_2$.   To complete the proof, the cocycle condition has to be verified for
the sheaves $\kf \otimes \kg$ and $\ke\!\xt^i(\kf,\kg)$.   It follows by
writing down the canonical isomorphisms arising by pulling back via the flat
morphisms $\sigma, \mu, p_2, p_{23}$ and their compositions.
\end{proof}
\vskip3mm

\begin{sub}\rm {\bf Invariance of supports}\label{invsup}\\
The annihilator ideal sheaf $\ka$ of any coherent $\ko_X$--module $\kf$ is
  the kernel of the homomorphism
\[
\ko_X \lra \kh\!\om_{\ko_X}(\kf, \kf)
\]
defined by $a \mapsto (s \mapsto as)$.   If $Y \xrightarrow{f} X$ is a flat
morphism, we obtain the exact diagram
\[
\UseComputerModernTips
\xymatrix{
0 \ar[r] & f^\ast \ka \ar[r]\ar@{-->}[d]^\approx & f^\ast \ko_X \ar[d]^\approx
\ar[r] & f^\ast \kh\!\om(\kf,\kf)\ar[d]^\approx\\
0 \ar[r] & \ka(f^\ast \kf) \ar[r] & \ko_Y \ar[r] & \kh\!\om(f^\ast \kf,
f^\ast \kf)
}
\]
with induced isomorphism $f^\ast \ka \xrightarrow{\approx} \ka (f^\ast \kf)$.

Applying this to $\sigma$ and $p_2$ of an action $G \times X
\xrightarrow{\sigma} X$ and an equivariant sheaf $\kf$, we obtain the
commutative diagram of isomorphisms
\[
\begin{array}{ccccc}

\sigma^\ast \ka & \cong &\ka(\sigma^\ast \kf) \cong \ka (p_2^\ast\kf)& \cong &

p_2^\ast \ka\\
\cap & & & &\cap\\
\sigma^\ast \ko_X & \cong &\ko_{G \times X} &\cong & p^\ast_2\ko_X\,.
\end{array}
\]
It is, again, straightforward to verify that the induced isomorphism
$\sigma^\ast\ka \cong p_2^\ast \ka$ satisfies the cocycle condition.
Consequently, the subscheme $\Supp (\kf) = \Supp(\ko_X/\ka)$ is
$G$--invariant.   In particular, $\sigma$ induces an action on $\Supp(\kf)$.
\end{sub}
\vskip3mm

\begin{sub}\rm{\bf Deficiencies of depth}\label{defs}\\ 
 Let $\kf$ be a coherent $\ko_X$--module on a scheme $X$, in our situation 
locally of finite type over $k$.   It is well--known that the sets
\[
S_m(\kf) = \{x \in X \mid \text{ prof } \kf_x \le m\}
\]
are closed.   This follows by locally embedding $X$ into an affine scheme
$\A^n$ and considering a finite free resolution of $\kf$ in $\A^n$.
Alternatively,
\[
S_m(\kf) = \underset{i \ge n-m}{\cup} \Supp\bigl(\ke\!\xt^i_{\ko_{X}} (\kf,
\ko_{X})\bigr)
\]
if $X$ is nonsingular of pure dimension $n$.   This can be used to define a
structure sheaf for $S_m(\kf)$.   In general, we consider $S_m(\kf)$ simply as
a reduced  subscheme.   By definition
\[
S_0(\kf) \subset S_1(\kf) \subset \ldots \subset S_m(\kf) \subset \ldots\,.
\]
It is well--known that $\dim S_m(\kf) \le m$ for any $m$, as in the analytic
case, see \cite{SiuTrm}.
\end{sub}
\vskip3mm

\begin{sub}{\bf Lemma:}\label{prof}
  Let $G \times X \xrightarrow{\sigma} X$ be an action of the linear algebraic
  group $G$ on $X$  and
  let $\kf$ be a $G$--equivariant coherent $\ko_X$--module.   Then every
  deficiency set $S_m(\kf)$ is $G$--invariant.
\end{sub}

\begin{proof}
  By Corollary \ref{invsub}, we only have to show that $S_m(\kf)$ is invariant
  under each isomorphism $X \xrightarrow{g} X$, $g \in G(k)$.   Because
  $g^\ast \kf \cong \kf$ for any $g$, we have, indeed,
\[
g^{-1} S_m(\kf) = S_m(g^\ast \kf) = S_m(\kf)\,.
\]
\end{proof}

\section{Sheaves of local cohomology}\label{locco}

Let $A \subset X$ be a reduced subscheme and let $\kf$ be coherent on $X$.
  We denote by $\kh^0_A \kf \subset \kf$ the subsheaf of all germs, which are
  annihilated by some power of the ideal sheaf $\ki_A$ of $A$.   Then
\[
\kh^0_A \kf = \underset{n}{\varinjlim}\;\; \kh\!\om(\ko/\ki^n_A, \kf)
\]
and $\kh^0_A\kf$ is coherent.  It is then straightforward to
prove that this equality extends to the derived functors, 
\[
\kh^i_A \kf \cong \underset{n}{\varinjlim}\;\; \ke\!\xt^i (\ko/\ki^nA_, \kf)\,,
\]
see \cite[\'{E}xpos\`{e}]{SGA2}.

In general, the sheaves $\kh^i_A\kf$ are only quasi--coherent, but not
coherent.   The following criterion on their vanishing and coherence can be
found in \cite{SGA2} and, for the analytic case, in \cite{SiuTrm}.
\vskip3mm

\begin{sub}{\bf Criterion:} \label{crit1}
  Let $\kf$ be a coherent $\ko_X$--module, $A \subset X$ a closed subset and
  $q \ge 0$ an integer.   Then
  \begin{enumerate}
  \item $\kh^i_A\kf = 0$ for $0 \le i \le q$ if and only if $\dim A \cap
  S_{m+q+1} (\kf) \le m$ for all $m$.
\item $\kh^i_A\kf$ is coherent for $0 \le i \le q$ if and only if $\dim A \cap
  \overline{S_{m+q+1}(\kf|X\smallsetminus A)} \le m$ for all $m$.
  \end{enumerate}
\end{sub}
\vskip3mm

\begin{sub}{\bf Lemma:}\label{locoh}
  Let $G \times X \xrightarrow{\sigma} X$ be an action as above, let $A
  \subset X$ be a closed invariant subset and let $\kf$ be a $G$--equivariant
  $\ko_X$--module.   Then the sheaves $\kh_A^i\kf$ are all $G$--equivariant.
\end{sub}
\begin{proof}
  The ideal sheaf $\ki_A$ of $A$ is equivariant by the assumption on $A$.
  The exact sequences
\[
\UseComputerModernTips
\xymatrix@R=8pt{
\ko \ar[r] & \ki^{n+1}_A \ar[r] & \ki^n_A \ar[r] & \ki^n_A/\ki_A^{n+1}
\ar@{=}[d] \ar[r] & 0\\
& & & \ki^n_A \otimes \ko_X/\ki_A & 
}
\]
may be used to prove that, also, all the powers $\ki^n_A$ are $G$--equivariant
by induction, see Lemma \ref{indequi}.  By the same lemma, also, the sheaves
$\ke\!\xt^i_{\ko_X} (\ko_X/\ki^n_A,\kf)$ are $G$--equivariant.  Finally, we
have the commutative diagrams
\[
\UseComputerModernTips
\xymatrix{
\sigma^\ast \ke\!\xt^i(\ko_X/\ki^n_A,\kf) \ar[r] \ar[d]^{\Phi_n}_\approx &
\sigma^\ast \ke\!\xt^i(\ko_X/\ki^{n+1}_A,\kf)\ar[d]_\approx^{\Phi_{n+1}}\\
p_2^\ast \ke\!\xt^i(\ko_X/\ki^n_A,\kf) \ar[r] & p_2^\ast \ke\!\xt^i(\ko_X/
\ki^{n+1}_A,\kf)
}
\]
where the isomorphisms $\Phi_n$ are induced by the corresponding compatible
isomorphisms for the sheaves $\ko_X/\ki^n_A$.   By this, we obtain an induced
isomorphism
\[
\sigma^\ast \kh^i_A \kf \underset{\approx}{\lra} p_2^\ast \kh^i_A \kf
\]
on the inductive limits.   A canonical, but elaborate, diagram yields the
cocycle condition for these isomorphisms.
\end{proof}
\vskip3mm

\begin{sub}\rm {\bf Subsheaves of bounded torsion}\label{btor}\\  
For any coherent module $\kf$ on $X$ and any dimension $d$, we define the
subsheaves $\kh^0_d \kf \subset \kf$ as
\[
\kh^0_d\kf = \varinjlim \kh^0_A\kf\,,
\]
the limit to be taken over all closed subsets $A \subset X$ of dimension $\le
d$.   Because $\dim S_d(\kf) \le d$ and $(\kh^0_A\kf)_x = 0$ for $x \in A
\smallsetminus S_d(\kf)$, $\dim_x A \le d$, the limit is stationary and
\[
\kh^0_d\kf = \kh^0_{S_d(\kf)} \kf\,.
\]
Therefore, the sheaves $\kh^0_d\kf$ are coherent for any $d$, and
$G$--invariant if $\kf$ is $G$--equivariant.  
\end{sub}
\vskip3mm

\begin{sub}\rm {\bf Remark:}\label{dertor}
There are also the sheaves $\kh^i_d\kf$, defined by the
derived functors $\kh^i_d$ of $\kh^0_d$. 
There is an analogous criterion for vanishing and coherence of these sheaves.
Then one can prove that for a coherent and $G$-equivariant $\kf$ also
$\kh^i_d\kf$ is $G$--equivariant if it is coherent.
\end{sub}

\section{Equivariant primary decomposition}\label{sectepd}
In this section, $X$ will be a scheme, locally of finite type over $k$, $G
\times X \xrightarrow{\sigma} X$ an action of a linear algebraic group over
$k$ and $\km$ a coherent $\ko_X$--module.   The main result is Theorem
\ref{epd} on the equivariant primary decomposition, for which we assume that
$G$ is connected.   The more general result for non--connected $G$ will be
derived from the former in Theorem \ref{ncepd}.
\vskip3mm

\begin{sub}\rm{\bf Invariance of irreducible components}\label{invirr}\\ 
When $G$ is connected, any irreducible component $Y_0$ of a $G$--invariant
closed subscheme $Y$ is again $G$--invariant:  let $y \in Y_0(k)$ and let $A =
\overline{o(y)}$ be the closure of the orbit.   Let $Y = Y_0 \cup Y_1$ and let
$A_i = A \cap Y_i$.   Because $A \subset Y$, we have $A = A_0 \cup A_1$ with
$A_0 \not= \emptyset$, and because $A$ is irreducible, we conclude that $A =
A_0 \subset Y_0$.   It follows that $g \cdot Y_0 = Y_0$ for any $g \in G(k)$.
\par
When $G$ is not connected, $G$ is a finite union $G = G^0 \cup h_1G^0 \cup
\dots \cup h_n G^0$, where $G^0$ is the identity component of $G$.   Let $H =
\{e, h_1, \dots, h_n\}$.   If $Y_0 \subset Y$ is an irreducible component of
the $G$--invariant closed subscheme $Y$, we consider the $G$--invariant
closure
\[
H. Y_0 = \underset{h \in H}{\cup} h. Y_0\,.
\]
Because $Y_0$ is $G^0$--invariant, $H.Y_0$ is the smallest $G$--invariant
closed subset containing $Y_0$.   It is then a finite union of irreducible
$G^0$--invariant components of the same dimension.   We call it a
$G$--component of $Y$.
\end{sub}
\vskip3mm

\begin{sub}\rm {\bf Associated components}\label{asscom}\\ 
We let $\Ass(\km)$ be the set of all closed irreducible subsets, which occur
as a component of $\Supp \kh^0_d(\km)$ for some $d$.   Because each
$\kh^0_d\km$ is coherent, the system $\Ass(\km)$ is locally finite.
Moreover, if $\km$ is $G$--equivariant, any $Y \in \Ass(\km)$ is
$G$--invariant, because also each $\kh^0_d\km$ is then $G$--equivariant and so
its support is $G$--invariant.
\end{sub}
\vskip3mm

\begin{sub}{\bf Remark:}\rm \label{pcyc}\; 
  The system $\Ass(\km)$ is the system of primary cycles of $\km$ in the sense
  of \cite[3.1.1]{EGAIV}.   There the set $\Ass_0(\km)$ is defined as the
  set of all points $x \in X$, closed or not, for which the maximal ideal
  $\fm_x \subset \ko_{X,x}$ belongs to $\Ass(\km_x)$ of the module $\km_x$.
  An irreducible closed subset $Y = \overline{\{x\}}$ is called a primary
  cycle if $x \in \Ass_0(\km)$.   If, by abuse of notion, we identify 
$\Ass_0(\km)$
  with the set of primary cycles, we get
\[
\Ass(\km) = \Ass_0(\km)\,.
\]
For a proof, one can use the criterion in \cite[3.1.3]{EGAIV}, which says that
$Y = \overline{\{x\}}$ is a primary cycle if and only if there is an open
neighbourhood $U(x)$ and a coherent subsheaf $\kg \subset \km|U$ such that $Y
\cap U$ is a component of $\Supp(\kg)$, or equivalently, if and only if there
is an open neighbourhood $U(x)$ and a section $s \in \Gamma(U,\km)$, such that
$Y \cap U$ is a component of $\Supp(\ko_Us)$.   If, now, $Y \in \Ass(\km)$,
then $Y \in \Ass_0(\km)$ using $\kg = \kh^0_d\km$.   Conversely, if $Y \cap U$
is a component of $\supp(\ko_Us)$ and $d = \dim Y$, then
\[
Y \cap U \subset U \cap \Supp \kh^0_d(\km)\,,
\]
and because $\dim\Supp\kh^0_d(\km) \le d$, it follows that $Y$ is a component
of $\supp\kh^0_d(\km)$.
\end{sub}
\vskip3mm

\begin{sub}{\bf Primary submodules}\label{primes}\\ \rm 
Let $\kn \subset \km$ be a coherent submodule.
$\kn$ is called $Y$--primary if $\Ass(\km/\kn) = \{Y\}$.   By definition of
$\Ass$ we obtain immediately

\begin{subsub}{\bf  Lemma:}\label{prim1}\textit{ $\kn$ is $Y$--primary if and 
only
if $Y = \Supp(\km/\kn)$ is irreducible and $\kh^0_d(\km/\kn) = 0$ for $d <
\dim Y$.}
\end{subsub}

\begin{subsub}{\bf Lemma:}\label{ass2}
\textit{ $\Ass(\km) \not= \emptyset$ if and only if $\km \not= 0$.}
\end{subsub}
\end{sub}
\vskip3mm

\begin{sub}{\bf Remark}\label{assrem}\rm \; Let $A$ be a noetherian ring, 
$N \subset M$
  finitely generated $A$--modules, $\fp \subset A$ a prime ideal.  Then the
  subsheaf $\widetilde{N}$ of $\widetilde{M}$ on $\Spec(A)$ is
  $V(\fp)$--primary, if and only if $N \subset M$ is $\fp$--primary in the
  sense of commutative algebra.  This follows by standard arguments
  considering the $A$--modules corresponding to
  $\kh^0_d(\widetilde{M}/\widetilde{N})$.
\end{sub}
\vskip3mm

\begin{sub}{\bf Lemma:}\label{prim2}
With notation as above, the following conditions are equivalent:
\begin{enumerate}
  \item[(a)] $\kn \subset \km$ is $Y$--primary.
  \item[(b)] For any open subset $U \subset X$ with $U \cap Y \not= \emptyset$,
  $\kn|U \subset \km|U$ is $U \cap Y$--primary.
  \item[(c)] $\kn(U) \subset \km(U)$ is $\ki_Y(U)$--primary for any affine open
  subset $U \subset X$.
  \end{enumerate}
\end{sub}
\begin{proof}
  (a) and (b) are equivalent by definition.   (c) follows from (b) by Remark
      \ref{assrem}.   (c) implies (a) by Lemma \ref{prim1} because
      $\Supp(\km/\kn) = Y$ and $\kh^0_d(\km/\kn) = 0$ for $d < \dim Y$ follow
      from the same local statements.
\end{proof}
\vskip3mm

\begin{sub}{\bf Main Lemma.}\label{mainlem} \textit{Let $Y \in \Ass(\km)$
    be a primary component of $\km$.  Then there is a coherent submodule $\kn
    \subset \km$ with $\Ass(\kn) = \{Y\}$.  If $G$ is connected and if $\km$
    is $G$--equivariant, then $\kn$ can be chosen to be $G$--equivariant.}
\end{sub}
\begin{proof} (a)\; Let $d = \dim Y$.   Then $Y$ is a component of $\Supp
  \kh_d^0\km$.   To see this, let $\kappa$ be minimal such that $Y$ is a
  component of $\Supp \kh^0_\kappa \km$.   There is an open subset $U$ such
  that $U \cap Y = U \cap \Supp \kh^0_\kappa \km$, and then $U \cap Y = U \cap
  \Supp \kh^0_d\km$ by definition of $\kh^0_\kappa \km$ as subsheaf of torsion
  of $\dim \le \kappa$.   We also have
\[
0 \not= \kh^0_Y \km \subset \kh^0_d\km\,.
\]
Let $\kf = \kh^0_Y \km$, let $\kt = \kh^0_{d-1} \kf$ and $A = \supp \kt$.
Let, further,
\[
\widetilde{A} = \underset{i \in I}{\cup} A_i
\]
be the union of the $(d-1)$--dimensional irreducible components of $A$.   This
system is locally finite and $\dim A_i = d-1$ for any $i$ because $\dim \Supp
\kh^0_{d-1} \kf \le d-1$.   Let $\eta \in Y$ and $\alpha_i \in A_i$ denote the
generic points, $\eta \not\in A$.   There is an affine neighbourhood $U_i$ of
$\alpha_i$ and an integer $m_i$ such that
\[
\ki^{m_i}_{A_i} \kt = 0 \text{\; on\; } U_i\,.
\]
By the Lemma of Artin--Rees, applied to the affine scheme $U_i$, there is an
integer $q_i$ such that
\[
(\ki^n_{A_i} \kf) \cap \kt = \ki^{n-q_i}_{A_i} \bigl((\ki^{q_i}_{A_i} \kf)
\cap \kt \bigr) \subset \ki^{n-q_i}_{A_i} \kt
\]
over $U_i$.   Hence there are integers $n_i > q_i$ such that
\[
(\ki^{n_i}_{A_i} \kf) \cap \kt = 0 \text{\; on\; } U_i\,.
\]
Because the system $(A_i)_{i\in I}$ is locally finite, the sheaf
\[
\kg = \underset{i \in I}{\cap} \ki^{n_i}_{A_i} \kf
\]
is coherent and $\kg \cap \kt = 0$ along $U = \underset{i}{\cup}\, U_i$.
Moreover, $\kg | Y \smallsetminus \widetilde{A} = \kf$, hence $\Supp \kg =
\Supp \kf = Y$.  Let $A^\prime = A \smallsetminus U$.  Because all the generic
points $\alpha_i \in U$, we have $\dim A^\prime \le d-2$.  Moreover,
$\kh^0_{d-1} \kg = 0$ on $U$ and then over $Y \smallsetminus A^\prime$.
\par
(b)\; Applying the same procedure to $\kg$ and $A^\prime$, etc., we can do a
finite induction to arrive at a coherent subsheaf $\kn \subset \kf$ with
$\Supp(\kn) = Y$ and $\kh^0_{d-1} \kn = 0$ on $Y$.   This implies that
\[
\Ass(\kn) = \{Y\}\,.
\]
\par
(c)\; If $\km$ is $G$--equivariant, then so is $\kh^0_Y \km = \kf$ above and
$\kh^0_{d-1}\kf$.   Then all the components $A_i$ are $G$--invariant.   It
follows that each $\ki^{n_i}_{A_i} \kf$ is $G$--equivariant, too, and that,
finally, $\kg$ is $G$--equivariant.   This shows that $G$--equivariance can be
maintained in the induction process (b).   
\end{proof}
\vskip3mm

The following two lemmata can be proven by the same standard arguments
as in commutative algebra for the analogous statements.

\begin{sub}{\bf Lemma:}\label{quass}
 Let $\kn \subset \km$ be a coherent submodule.  Then $\Ass(\kn) \subset
  \Ass(\km) \subset \Ass(\kn)$ $\cup \Ass(\km/\kn)$.
\end{sub}
\vskip3mm

\begin{sub}{\bf Corollary:}\label{mass}
  Let $\kn_1, \kn_2 \subset \km$ be coherent submodules, and suppose that
  there are subsets $S_1, S_2 \subset \Ass(\km)$ such that $\Ass(\kn_i) =
  \Ass(\km)\smallsetminus S_i$ and $\Ass(\km/\kn_i) = S_i$ for $i = 1,2$.
  Then
\[
\Ass(\kn_1 \cap \kn_2) = \Ass(\km) \smallsetminus S_1 \cup S_2 \text{ and }
\Ass(\km/\kn_1 \cap \kn_2) = S_1 \cup S_2\,.
\]
\end{sub}
\vskip3mm

\begin{sub}{\bf Proposition:}\label{arbass}
  Let $S \subset \Ass(\km)$ be any subset.   Then there is a coherent
  submodule $\kn \subset \km$ such that
\[
\Ass(\kn) = \Ass(\km) \smallsetminus S \text{\; and\; } \Ass(\km/\kn) = S\,.
\]
If $G$ is connected and if $\km$ is $G$--equivariant, $\kn$ can be chosen as a
$G$--equivariant submodule.
\end{sub}

\begin{proof}
  We use the Lemma of Zorn, applied to the following set $\E$ of coherent
  submodules.
\[
\E := \{\kn \subset \km \mid \kn \text{ coherent } G\text{--invariant
  submodule s.t. } \Ass(\kn) \subset \Ass (\km) \smallsetminus S\}\,.
\]
$\E \not= \emptyset$ because $\kn = 0$ belongs to $\E$.   The set $\E$ is
inductively ordered, i.e.\ for any chain
\[
\kn_1 \subset \kn_2 \subset \ldots
\]
of members of $\E$ there is a supremum in $\E$.   This is the union $\kn =
\underset{\nu}{\cup} \kn_\nu$.   Indeed, $\kn$ is coherent because $X$ is
locally noetherian, and the chain terminates locally.   Because $\sigma^\ast
\km \xrightarrow{\Phi} p_2^\ast \km$ maps each $\sigma^\ast \kn_\nu$
isomorphically to $p_2^\ast \kn_\nu$, also $\kn$ is a $G$--equivariant
submodule.   Moreover,
\[
\Ass(\kn) = \underset{\nu}{\cup} \Ass(\kn_\nu)\,.
\]
To see this, let $Y \in \Ass(\kn)$, $\dim Y = d$, $Y$ a component of $\supp
\kh^0_d \kn$.   We have
\[
\kh^0_d\kn_\nu \subset \kh^0_d \kn\,,
\]
all of dimension $\le d$.  Let $\eta \in Y$ be the generic point.   Because
the sequence terminates at $\eta$, there is an index $\nu$ with $\eta \in
\Supp \kh^0_d \kn_\nu$.   Then $Y \in \Ass(\kn_\nu)$.

Since $\E$ is now inductively ordered, there is a maximal $\kn \in \E$.  We
verify that then $\Ass(\km/\kn) \subset S$.  Let $Y \in \Ass(\km/\kn)$.  By
the Main Lemma there is a submodule $\kf \subset \km$, containing $\kn$ such
that $\Ass (\kf/\kn) = \{Y\}$.  Then $\kf \not\in \E$, i.e.\ $\Ass(\kf)
\not\subset \Ass(\km)\smallsetminus S$.  On the other hand, 
$\Ass(\kf) \subset \Ass(\kn)\cup \Ass(\kf/\kn) = \Ass(\kn) \cup \{Y\}$.  
This proves that $Y \in S$.
Let, conversely, $Y \in S$.   Then $Y \not\in \Ass(\kn)$, hence $Y \in
\Ass(\km/\kn)$.   We have shown that
\[
\Ass(\km/\kn) = S\,.
\]
Finally, let $Y \in Ass(\km) \smallsetminus S$.   Then $Y \not\in
\Ass(\km/\kn)$ and, hence, $Y \in \Ass(\kn)$.   This proves that
\[
\Ass(\kn) = \Ass(\km) \smallsetminus S\,.\qed
\]
\renewcommand{\qed}{}\end{proof}
\vskip3mm

\begin{sub}{\bf Theorem:} ({\rm Equivariant Primary Decomposition})\label{epd}
  Let $X$ be a scheme, locally of finite type over $k$, let $G \times X
  \xrightarrow{\sigma} X$ be an action of a connected linear algebraic group
  over $k$ on $X$ and let $\km$ be a coherent $G$--equivariant
  $\ko_X$--module.   Then
  \begin{enumerate}
  \item For any $G$--equivariant coherent submodule $\kn \subset \km$, the
  system $\Ass(\km/\kn) = \{Y_i\}_{i \in I}$ is locally finite, consisting of
  closed $G$--invariant subsets $Y_i$.
\item There is a family $(\kq_i)_{i \in I}$ of $G$--equivariant coherent
  modules, $\kn \subset \kq_i \subset \km$, such that each $\kq_i$ is
  $Y_i$--primary and such that
\[
\kn = \underset{i \in I}{\cap} \kq_i\,,
\]
this intersection being locally finite.
  \end{enumerate}
\end{sub}
\begin{proof}
  By Proposition \ref{arbass} there is a $G$--equivariant coherent submodule
  $\kq_i \subset \km$ such that $\Ass(\km/\kq_i) = \Ass
  \bigl((\km/\kn)/(\kq_i/\kn)\bigr) = \{Y_i\}$ and $\Ass(\kq_i/\kn) =
  \Ass(\km/\kn) \smallsetminus \{Y_i\}$.   Then $\kq_i$ is $Y_i$--primary.
  Because $\Supp \km/\kq_i = Y_i$, the intersection
\[
\kn^\prime = \underset{i \in }{\cap} \kq_i
\]
is locally finite, coherent and $G$--equivariant.   We have
\[
\kn^\prime/\kn = \underset{i \in I}{\cap} (\kq_i/\kn)\,.
\]
Because $\Ass(\kn^\prime/\kn) \subset \Ass(\km/\kn) \smallsetminus \{Y_i\}$
for any $i \in I$, we find $\Ass(\kn^\prime/\kn) = \emptyset$ and, hence, $\kn
= \kn^\prime$.                                    
\end{proof}
\vskip3mm

\begin{sub}{\bf Relative gap sheaves}\label{gap}\\ \rm 
The primary submodule $\kq_i$ is called unembedded, if the primary component
  $Y_i$ is not contained in another primary component. In Corollary \ref{unemb}
  we will prove that the unembedded primary submodules are unique.
  For this, we introduce the so--called gap sheaves $\kn[A] \subset \km$ of a
  coherent submodule $\kn \subset \km$ with respect to a closed subset $A
  \subset X$ as the submodule with
\[
\kn[A]/\kn = \kh^0_A(\km/\kn)\,.
\]
Similarly, the submodule $\kn[d] \subset \km$ is defined as the submodule with
\[
\kn[d]/\kn = \kh^0_d(\km/\kn)\,,
\]
see \cite{SiuTrm}.   These gap sheaves are both coherent.   If $\kn$ is a
$G$--equivariant submodule of the $G$--equivariant module $\km$, then
$\kn[d]$ is $G$--equivariant as well as $\kn[A]$ if $A$ is $G$--invariant.
The following proposition has been proved by Y.T.~Siu in \cite{Siu} in the
analytic case.  
\end{sub}
\vskip3mm

\begin{sub}{\bf Proposition:}\label{gapdec}
  Let $\kn = \cap \kq_i$ be a primary decomposition of the coherent submodule
  $\kn$ of $\km$ as in Theorem \ref{epd}.   Then:
  \begin{enumerate}
  \item For any closed subset $A \subset X$ the gap sheaf $\kn[A]$ has the
  decomposition
\[
\kn[A] = \underset{Y_i \not\subset A}{\cap} \kq_i\,.
\]
\item For any integer $d \ge 0$ the gap sheaf $\kn[d]$ has the decomposition
\[
\kn[d] = \underset{d < \dim Y_i}{\cap} \kq_i\,.
\]
\end{enumerate}
\end{sub}
\begin{proof}
  It follows easily from the definition of $\kn[A]$ that $\kn[A] = \cap
  \kq_i[A]$, using the local finiteness of the system $\{\kq_i\}$.   For a
  fixed index $j$, we have $\Supp(\km/\kq_j) = Y_j$ and $\kh^0_d(\km/\kq_j) =
  0$ for $d < \dim Y_j$.   Therefore, $\kh^0_A (\km/\kq_j) = 0$ if $Y_j
  \not\subset A$ and $\kh^0_A(\km/\kq_j) = \km/\kq_j$ if $Y_j \subset A$.   It
  follows that
\[
\kq_j[A] = 
\begin{cases}
  \kq_j & \text{if $Y_j \not\subset A$}\\
  \km   & \text{if $Y_j \subset A$}\,.
\end{cases}
\]
This proves (1).   The statement (2) follows by the same argument because
\[
\kq_j[d] = 
\begin{cases}
  \kq_j & \text{if $d < \dim Y_j$}\\
  \km   & \text{if $\dim Y_j \le d$}\,.
\end{cases}
\]
\end{proof}
\vskip3mm

\begin{sub}{\bf Corollary:}\label{unemb}
  Let $\kn = \cap \kq_i$ be a primary decomposition, as in Theorem \ref{epd},
  and let $Z_j = \underset{i \not= j}{\cup} Y_i$.  If the primary component
  $Y_j$ is unembedded, then $\kq_j$ is unique and equal to $\kn[Z_j]$.
\end{sub}
\begin{proof}
  By \ref{gapdec} $\kn[Z_j] = \underset{Y_i \not\subset Z_j}{\cap} \kq_i$ and
  the intersection reduces to $\kq_j$ if $Y_j$ is unembedded.
\end{proof}
\vskip3mm

\begin{sub}{\bf Remark:}\label{noeth}\; \rm 
  If $X$ is of finite type over $k$, then $\Ass(\km/\kn)$ is finite because
  $X$ is noetherian.   In this case, there are convenient composition series,
  as in the following proposition, see \cite{EGAIV}, \cite{NBC}.
\end{sub}
\vskip3mm

\begin{sub}{\bf Proposition:}\label{cser}
  Let $\Ass(\km) = \{Y_1, \dots, Y_n\}$ be finite with any order of the
  primary components.   Then there is a (composition) series
\[
0 = \km_0 \subset \km_1 \subset \dots \subset \km_n = \km
\]
of coherent submodules such that $\km_{\nu-1}$ is primary in $\km_\nu$ with
$\Ass(\km_\nu/\km_{\nu-1}) = \{Y_\nu\}$ for $\nu = 1, \dots, n$.  If $G$ is
connected and if $\km$ is $G$--equivariant, the $\km_\nu$ can be chosen to be
$G$--equivariant.
\end{sub}
\begin{proof}
  By Proposition \ref{arbass},  there is a $Y_n$--primary submodule 
$\km_{n-1} \subset \km$
  such that $\Ass (\km/\km_{n-1}) = Y_n$ and $\Ass(\km_{n-1}) = \{Y_1, \dots,
  Y_{n-1}\}$.   Now proceed by descending induction.   In each step, $\km_\nu$
  can be chosen $G$--equivariant.
\end{proof}
\vskip3mm

\begin{sub}\rm {\bf Associated $G$--components}\label{assgc}\\
We now consider the case of an arbitrary linear algebraic group $G$ and a
$G$--equivariant coherent module $\km$ on $X$.   

Let $G^0$ be the identity component and let $G = G^0 \cup h_1 G^0 \cup \dots
\cup h_m G^0$ such that the residue classes $\bar{e}$, $\bar{h}_1, \dots
\bar{h}_m$ are the elements of the finite group $H = G/G^0$.   If $Y_0 \in
\Ass(\km)$, then $Y_0$ is $G^0$--invariant and we are given the
$G$--component
\[
Y = Y_0 \cup h_1 Y_0 \cup \dots \cup h_m Y_0\,.
\]
If $Y_0$ is a $d$--dimensional component of $\Supp \kh^0_d \km$, then $Y
\subset \Supp \kh^0_d\km$ and $Y$ is of pure dimension $d$.   We call such a
union an \textbf{associated $\mathbf{G}$--component} of $\km$.   A
$G$--equivariant coherent submodule $\kq \subset \km$ is called 
$\mathbf{Y}$\textbf{--primary},
if $\Ass(\km/\kq) = \{Y_1, \dots, Y_n\}$, where $Y_1, \dots, Y_n$ are the
irreducible components of the $G$--component of $Y$.
\end{sub}
\vskip3mm

\begin{sub}{\bf Theorem:}\; ({\rm Equivariant primary decomposition for 
arbitary $G$})\label{ncepd}\\
Let the linear algebraic group $G$ act on $X$ as stated in our general
assumption, and let $\km$ be a coherent $G$--equivariant $\ko_X$--module,
  let $\kn \subset \km$ be a coherent $G$--equivariant submodule, and let
  $(Z_i)_{i \in I}$ be the (locally finite) family of associated
  $G$--components of $\km/\kn$.   There is a family $(\kq_i)_{i \in I}$ of
  coherent $G$--equivariant submodules, $\kn \subset \kq_i \subset \km$, such
  that each $\kq_i$ is $Z_i$--primary and such that
\[
\kn = \underset{i \in I}{\cap} \kq_i\,,
\]
this intersection being locally finite.
\end{sub}
\begin{proof}
  We may assume that $\kn = 0$, and we let $\Ass(\km) = \{Y_j\}_{j\in J}$ with
  each $Y_j$ irreducible and $G^0$--invariant.   As in the proof of Theorem
  \ref{epd}, there are $Y_j$--primary $G^0$--equivariant submodules $\kq_j$
  satisfying $\Ass(\kq_j) = \Ass(\km) \smallsetminus \{Y_j\}$, and which
  constitute the $G^0$--equivariant primary decomposition of $0$ in $\km$.
  Now let $Z = Z_i$ be one of the associated $G$--components of $\km$, being a
  finite equidimensional union
\[
Z = Y_{j_1} \cup \dots \cup Y_{j_n}\,.
\]
We consider the intersection
\[
\kq = \kq_{j_1} \cap \dots \cap \kq_{j_n}
\]
which is $G^0$--equivariant and satisfies
\[
\Ass(\kq) = \Ass(\km) \smallsetminus S_Z \text{ and } \Ass(\km/\kq) = S_Z
\]
where $S_Z = \{Y_{j_1}, \dots, Y_{j_n}\}$, see Corollary \ref{mass}.
\par
We are going to replace $\kq$ with a $G$--equivariant submodule $\kq_Z$ with
$\Ass(\kq_Z) = \Ass(\kq)$ and $\Ass(\km/\kq_Z) = \Ass(\km/\kq)$.   For that,
let $h_1, \dots, h_m \in G$ be as above, and define the submodule $\kq(h_\mu)
\subset \km$ as the image of $h^\ast_\mu\kq$ under the isomorphism 
$\Phi_{h_\mu}$ with diagram
\[
\UseComputerModernTips
\xymatrix@1@=0pt@!{
h_\mu^\ast \kq \ar[d]_\approx & \subset &  h^\ast_\mu \km
\ar[d]^{\Phi_{h_\mu}}_{\approx}\\
\kq(h_\mu)& \subset & \km\,. 
}
\]
Now we put
\[
\kq_Z := \kq \cap \kq (h_1) \cap \dots \cap \kq(h_m)\,.
\]
By definition, $\Ass\bigl(\kq(h_\mu)\bigr) = \Ass(\km) \smallsetminus S_Z$ and
$\Ass\bigl(\km/\kq(h_\mu)\bigr) = S_Z$ for any $\mu$, and then the same holds
for $\kq_Z$ by Corollary \ref{mass}.   We are going to verify that $\kq_Z$ is
$G$--equivariant.   This proves the theorem because the sets $S_Z$ constitute
a (disjoint) partition of $\Ass(\km)$ as $Z$ varies in the family 
$(Z_i)_{i\in I}$. 

In order to show that $\kq_Z$ is $G$--equivariant, it is sufficient to prove
that for any $h \in \{h_1, \dots, h_m\}$ the isomorphism $\Phi_h$ maps
$h^\ast\kq_Z$ to $\kq_Z$, see Proposition \ref{inv1}.   By the cocycle
condition for $\Phi$, we have the diagrams
\[
\UseComputerModernTips
\xymatrix{
(h_\mu h)^\ast \kq \ar@/_2pc/[rrr]_{\Phi_{h_{\mu}h}}
&\!\!\!\!\!\!\!\!\!\!\!\!\! \cong  h^\ast h_\mu^\ast \kq 
 \ar[r]^>>>>>{h^\ast \Phi_{h_\mu}} &
h^\ast \kq(h_\mu) \ar[r]^{\Phi_h} & \kq(h_\mu h)
}
\]
for any $\mu$.   Let $\bar{h}_\nu = \bar{h}_\mu\bar{h}$ in the group $G/G^0$.
then $h_\mu h = g h_\nu$ for some $g \in G^0$.   Because $\Phi_g$ maps
$g^\ast\kq$ to $\kq$, we have the analogous diagram
\[
\UseComputerModernTips
\xymatrix{
(g h_\nu)^\ast \kq \ar@/_2pc/[rrr]_{\Phi_{gh_{\nu}}}
&\!\!\!\!\!\!\!\!\!\!\!\!\! \cong  h^\ast_\nu  g^\ast \kq 
 \ar[r]^{h^\ast_\nu \Phi_g} &
h^\ast_\nu \kq \ar[r]^{\Phi_{h_\nu}} & \kq(h_\nu),
}
\]
proving that $\kq(h_\mu h) = \kq(h_\nu)$ and that
\[
h^\ast \kq(h_\mu) \overset{\Phi_h}{\underset{\approx}{\longrightarrow}}
\kq(h_\nu)\,.
\]
Since multiplication with $\bar{h}$ is a bijection of the group $G/G^0$, we
conclude that
\[
h^\ast \kq_Z \overset{\Phi_h}{\underset{\approx}{\longrightarrow}} \kq_Z\,.
\]
This finishes the proof of the theorem.
\end{proof}

\textbf{Remark 1.}\;
  The $G$--equivariant decomposition in Theorem \ref{ncepd} automatically
  yields the $G^0$--equivariant decomposition if one inserts the
  decompositions of the $\kq_Z$, according to the decompositions $Z = Y_1 \cup
  \dots \cup Y_n$ into $G^0$--invariant irreducible components.

\textbf{Remark 2.}\;
  An associated $G$--component $Z = Y_1 \cup \dots \cup Y_n$ is unembedded if
  and only if each of its components $Y_\nu$ is unembedded.   In that case
  $\kq_Z$ is uniquely determined by the same proof as for Corollary
  \ref{unemb}. Then $\kq_Z$ can be chosen as $\kq$ in the above
  proof, because then $\kq$ is already $G$--invariant.
\vskip5mm

\section{Some toric sheaves on $\P_2$}\label{ptwo}

\begin{sub}\label{torlb}\rm {\bf Toric structures on line bundles on $\P_n$}\\
Here we use the standard notation for toric varieties and equivariant 
divisorial sheaves on them as is shortly explained in Section \ref{toricsh}.
The torus $(k^\ast)^{n+1}/k^\ast$ acts naturally on $\P_n(k)$ via
$\sigma(\langle t\rangle, \langle x\rangle)=\langle t_0x_0, \ldots, t_n
x_n\rangle,$
depending on the choice of the homogeneous coordinates. The invariant
Weil-divisors for this action are then the divisors 
$$D(a)= \Sigma\; a_\nu H_\nu,\quad a=(a_0,\ldots, a_n)\in \Z^{n+1},$$
where $ H_\nu$ is the hyperplane with equation $x_\nu$. As defined in 
\ref{equidiv}, the invertible sheaf $\ko(D(a))$ comes with a canonical
equivariance $\Phi_a.$ If $D(a)$ is effective, the canonical section
$$x^a=x_0^{a_0}\cdot\ldots\cdot x_n^{a_n} $$
of $\ko(D(a)),\; a_\nu\geq 0, $ is $T$-invariant. 
In this case the section $x^a$ spans the space of $T$-invariant sections,
$$\Gamma(\P_n, \ko(D(a)))^T=k\cdot x^a.$$

There is an obvious isomorphism
\[
\Hom(\P_n, \ko(D(a)), \ko(D(b)))^T\cong\Gamma(\P_n, \ko(D(b)-D(a)))^T
\]
between the 0- or 1-dimensional spaces of $T$--invariant homomorphisms and 
$T$--invariant sections. If $a\leq b$, we write $\ko(D(a))\xrightarrow{x^{b-a}}
\ko(D(b))$ for the homomorphism defined by the section $x^{b-a}$.
\end{sub}
\vskip4mm

\begin{sub}\rm {\bf Sheaves with cubic support on $\P_2$}\label{cubsupp}\\
Here we consider examples of torus equivariant sheaves of
dimension 1 on the projective plane, considered as a toric surface with the
natural action of $T=(k^\ast)^3/k^\ast$ as in \ref{torlb}. More generally, 
we first describe stable
coherent sheaves $\kf$ on $\P_2$ with Hilbert polynomial $\chi\kf(m)=3m+1$. It
has been shown in \cite{FrTr} that any such sheaf has a resolution
\[
0 \lra 2\ko(-2) \overset{A}{\lra} \ko(-1) \oplus \ko \lra \kf \lra 0,
\]
 where the homomorphism is identified with the matrix
\[
A = \begin{pmatrix} z_1 & q_1\\z_2 & q_2\end{pmatrix}
\]
with $z_1, z_2$ linearly independent linear forms and $q_1, q_2$ quadratic
forms. Then $\det(A) = z_1q_2 - z_2q_1 \not= 0$ defines a cubic curve $C$,
the support of $\kf$.  
Let $p$ be the common zero of $z_1, z_2$. Then the unique section of $\kf$
gives rise to the exact sequence
\[
0\to \ko_C\to \kf\to k_p\to 0
\]
such that $p$ is the zero locus of this section. Thus, $C$ and $p\in C$ are
invariants of the sheaf, which even determine it, see \cite{FrTr}. If $\kf$ is
torus--equivariant, then $C_{red}$ is a torus invariant divisor contained in
the union $H_0\cup H_1\cup H_2$ of the coordinate lines and $p$ must be a fixed
point. More precisely, we have the 
\begin{subsub}{\bf Lemma:}\label{exequi} \textit{Let $\kf$ be as above. Then 
the following are equivalent
\begin{itemize}
\item[(1)] $\kf$ is $T$--equivariant
\item[(2)] $C_{red}\subset H_0\cup H_1\cup H_2$ and $ p$ is a singular point 
of $C$ and a coordinate (or fixed) point of $\P_2$
\item[(3)] The matrix $A$ can be chosen as  
\[
\left(
\begin{array}{cc}
z_1 & w_1w_2\\
z_2 & 0
\end{array}
\right)
\]
where $z_1, z_2, w_1, w_2$ belong to the set $\{x_0, x_1, x_2\}$ of 
homogeneous coordinates.
\end{itemize}
In addition, if condition (3) holds, $\kf$ has the resolution
\[
0\to \ko(-D-L_1+L_2)\oplus \ko (-D)\xrightarrow{A} \ko(-D+L_2)\oplus
\ko(-L_1+L_2)\to \kf\to 0,
\]
where $L_i$ is the line with equation $z_i$ , $D_i$ the
line with equation $w_i$, and where $D=D_1+D_2$. The matrix $A$ is then 
invariant and under the canonical linearizations of the terms of the 
resolution and thus induces a linearization of $\kf$.}
\end{subsub}

\begin{proof} If $\kf$ is $T$--equivariant, $p$ and $C$ are invariant, and this
implies (2). Assume now that (2) is satisfied. The equivalence of the
representing matrices is defined by
\[
\begin{array}{lllcl}
\left(
\begin{array}{cc}
\alpha & \beta\\
\gamma & \delta
\end{array} 
\right)
&
\left(
\begin{array}{cc}
z_1 & q_1\\
z_2 & q_2
\end{array}
\right)
&
\left(
\begin{array}{cc}
\lambda & 0\\
z & \mu
\end{array} 
\right)
& = &
\left(
\begin{array}{cc}
z'_1 & q'_1\\
z'_2 & q'_2
\end{array}
\right)
\end{array}
\]
corresponding to the automorphism groups of the terms of the
resolution. Because $z_1(p)=z_2(p)=0$ and $p$ is a coordinate point, we may
assume that $z_1, z_2$ are coordinates in $\{x_0, x_1, x_2\}$.
Since $C_{red}$ is contained in the coordinate triangle, we have
\[
z_1q_2-z_2q_1=\lambda u_1u_2u_3,
\]
where the $u_i$ are coordinates again. One of the $u_i,\; u_1$ say, must vanish
at $p$. Therefore, we may assume that $z_2=u_1$. This implies that $q_2$ is
divisible by $z_2,\; q_2=z_2w$, with a linear form $w$. Again by an equivalence
operation we may even assume that $q_2=0$. Then the new representing matrix
has the form 
\[
A=\left(
\begin{array}{cc}
z_1 & w_1w_2\\
z_2 & 0
\end{array}
\right)
\]
with $z_1, z_2, w_1, w_2$ coordinates. This proves (3) from (2). Finally, the
shape of $A$ in case (3) implies (1) as described in the second part of the 
lemma.
\end{proof}

\textbf{Remarks:}\\
1) The resolution of type (3) is a special case of an equivariant
resolution of a toric sheaf on a toric variety, cf \cite{P4}. 

2) A general matrix with monomial entries which are invariant with respect 
to summands $\ko(D)$ need not be invariant.
It is so, however, if all of its minors are monomials, too, see 
\ref{matcd} and \cite{CD}. This is the case in the above example.
\vskip3mm

\begin{subsub}\rm {\bf The toric locus}\label{tloc}\\ 
In the above case there are only
finitely many toric sheaves in the moduli space $M_{\P_2} (3m+1)$ of stable
coherent sheaves with Hilbert polynomial $3m+1$ for a fixed torus
action. However, condition (2) tells us, that there is a torus action on $\P_2$
for which $\kf$ is equivariant as soon as $C_{red}$ decomposes into lines with
a singular point of $C$. Recalling from \cite{FrTr} that $M_{\P_2} (3m+1)$ is
isomorphic to the tautological cubic of the Hilbert scheme of cubics, we find
that the sheaves $\kf$ in $M_{\P_2} (3m+1)$ which admit a toric structure,
form a $6$--dimensional closed subvariety in the $10$--dimensional moduli 
space.
\end{subsub}
\end{sub}
\vskip4mm

\begin{sub}\rm {\bf The primary decomposition:} \label{cubprd}\\ 
Let $\kf$ be the $T$--equivariant sheaf on $\P_2$ with representing matrix
$$A = 
\begin{pmatrix}
  x_1 & x_0x_1\\x_2 & 0 
\end{pmatrix}$$
\end{sub}
as in the previous section. Then $C$ is the union of the coordinate lines
$L_0, L_1, L_2$ with equation $x_0,x_1,x_2$ respectively. 
Here $\kf$ is Cohen--Macaulay and locally free on $C
\smallsetminus \{p_0\}$, where $p_0 = \langle 1,0,0\rangle$, etc.   We have
$\kh_0^0\kf = 0$ such that $\Ass(\kf) = \{L_1, L_2, L_3\}$.   By the
decomposition theorem there are unique $T$--equivariant coherent primary
submodules $\kf_0, \kf_1, \kf_2$ with $\kf_0 \cap \kf_1 \cap \kf_2 = 0$ and
$\Ass(\kf/\kf_\nu) = \{L_\nu\}$.   Then 
$\kf_0 = \kh^0_{L_1 \cup L_2} \kf$, $\kf_1
= \kh^0_{L_0 \cup L_2} \kf$, $\kf_2 = \kh^0_{L_0 \cup L_1} \kf$ by Corollary
\ref{unemb}. 

The sheaves $\kf_\nu$ and their quotients $\kf/\kf_\nu$ can explicitly be
presented in the diagrams
\[
\UseComputerModernTips
\xymatrix{
& & 0\ar[d] & 0 \ar[d] & \\
0 \ar[r] & 2\ko(-2)\ar@{=}[d] \ar[r]^{A_\nu} & 2\ko(-1) \ar[d]^{B_\nu} \ar[r]
& \kf_\nu \ar[d] \ar[r] & 0\\
0 \ar[r] & 2\ko(-2) \ar[r]^A & \ko(-1) \oplus \ko \ar[d] \ar[r] & \kf \ar[d]
\ar[r] & 0\\
& & \ko_{L_\nu} \ar@{=}[r]\ar[d] & \ko_{L_\nu} \ar[d] & \\
& & 0 & 0  
}
\]
with
\[
\begin{array}{ll}
A_0 = 
\begin{pmatrix}
  x_1 & x_1\\x_2 & 0
\end{pmatrix}, & B_0 = 
\begin{pmatrix}
  1 & 0\\0 & x_0
\end{pmatrix}\\[4ex]
A_1 = 
\begin{pmatrix}
  x_1 & x_0\\x_2 & 0
\end{pmatrix}, & B_1 = 
\begin{pmatrix}
  1 & 0\\0 & x_2
\end{pmatrix}\\[4ex]
A_2 = 
\begin{pmatrix}
  x_1 & \phantom{-}0\\x_2 & -x_0
\end{pmatrix}, & B_2 = 
\begin{pmatrix}
  1 & x_0\\0 & x_2
\end{pmatrix}\,.
\end{array}
\]
\vskip2mm

While $\kf_1, \kf_2$ are locally free on $L_0 \cup L_2$, resp.$\ L_0 \cup
L_1$, the sheaf $\kf_0 \cong \ko_{L_1} \oplus \ko_{L_2}$ is not free at the
intersection point $p_0$.
\vskip10mm

\section{Sheaves on  varieties admitting a homogeneous coordinate ring}\label{sqp}

In this section we consider the following setting which generalizes the
setting of a coordinate ring for a variety $X$. Let $Z$ be an integral affine
scheme with coordinate ring $S$, let
\[
G\times Z\xrightarrow{\sigma} Z
\]
be an action of a linear algebraic group on $Z$, let $H\subset G$ be a
diagonalizable closed normal subgroup, let $W\subset Z$ be a 
$G$--invariant open
subset and assume that there is a good quotient $X=W//H$ of $W$ by $H$. 
Note that the rings $S$ and $k[G]$ are then $X(H)$--graded via the dual 
actions of $H\times Z\to Z$ and $H\times G\to G.$, see \ref{dact}.
For later use we insert the 
\vskip5mm

\begin{sub}{\bf Lemma:}\label{lemdesc}
Let $M\overset{\vartheta}\lra k[G]\otimes M$ be a dual action
and let $M\overset{\vartheta_H}\lra k[H]\otimes M$ be the induced dual action 
of $H$. Then $\vartheta$ maps the component $M_\chi$ to 
$k[G]_\chi\otimes M$ for any character $\chi$ of $H.$  
\end{sub}
\vskip5mm

\begin{proof}

Let $\vartheta(m)=\sum \varphi_i\otimes m_i$ such that the
$\varphi_i$ are homogeneous with characters $\chi_i$.
Because $\varphi_i|H=c_i\chi_i$ with scalar factors $c_i=\varphi_i(e)$,
we obtain 
$$\chi \otimes m= \vartheta_H(m)=\Sigma_{i\in I_1} \chi_i\otimes c_im_i $$
where the index set $I_1$ is defined by $c_i\not =0$. 

Because the set $X(H)$ of characters of $H$ is a basis of the vector 
space 
$k[H]$, it follows that $\chi=\chi_i$ for all $i\in I_1$, assuming that all 
$m_i\not = 0,$ and we have then $m= \sum_{i\in I_1} c_i m_i$.

In order to get information about the $\varphi_i$ for 
$i\in I_0=I\smallsetminus I_1$, we 
use the 
rule $(1\otimes\vartheta)\circ\vartheta=(\mu^\ast\otimes 1)\circ\vartheta$
and the fact that $\mu^\ast(\varphi_i)(h,g)=(\chi_i\otimes\varphi_i)(h,g).$
Applying this to $\vartheta(m)$ and restricting to $k[H]\otimes k[G]\otimes M$,
we obtain the equalities
$$\chi\otimes\vartheta(m)=\chi\otimes\Sigma_{i\in I_1}\vartheta(c_im_i)=\chi\otimes\Sigma_{i\in I_1}\varphi_i\otimes m_i+\Sigma_{i\in I_0}\chi_i\otimes\varphi_i\otimes m_i.$$
Now we can shift terms of index $i\in I_0$ and $\chi_i=\chi$ to 
the first sum and finally obtain 
$$\vartheta(m)=\Sigma_{j\in J}\varphi_j\otimes m_j$$
with all $\varphi_j\in k[G]_\chi.$
\end{proof}

We are going to associate to any $X(H)$--graded $S$--module $M$ a
quasicoherent sheaf $\widetilde{M}$ on $X$ and to any dual action
$M\xrightarrow{\vartheta} k[G]\otimes_k M$ which induces the dual action 
of $H$ on $M$ corresponding to the $X(H)$--grading, a $G$--linearization of
$\widetilde{M}$ with the properties listed in Proposition \ref{shff}.

\begin{sub}{\bf Construction of $\mathbf{\widetilde{M}}$:}\label{shm}\rm

Let $\vartheta_H$ be the dual $H$--action given by the $X(H)$-grading and let
\[
\sigma^\ast_H \km\xrightarrow{\varphi_H} p^\ast_{2H} \km
\]
be the corresponding $H$--linearization of the sheaf $\km$ on $Z$,
where the index $H$ refers to the induced action 
$H\times Z\xrightarrow{\sigma_H} Z$, see \ref{dactdiag}.

For any open $U\subset X$ let $\widetilde{U}\subset W$ denote its preimage
under the quotient map $W\xrightarrow{\pi} X$. The restriction of
$\varphi_H$ to $H\times \widetilde{U}$ defines a dual action
\[
\km(\widetilde{U})\to k[H]\otimes_k\km(\widetilde{U})
\]
which in turn defines an $X(H)$--grading
\[
\km(\widetilde{U})=\underset{\chi}{\oplus}\km(\widetilde{U})_\chi\ .
\]
These gradings are compatible with the restrictions to open subsets and thus
define an $X(H)$--grading of the sheaf
\[
\pi_\ast(\km|W)=\underset{\chi}{\oplus}\widetilde{\km}_\chi
\]
where $\widetilde{\km}_\chi(U)=\km(\widetilde{U})_\chi$. Let
$\widetilde{M}=\widetilde{\km}_0$ be the degree zero part of that sheaf. It is
called the sheafification of $M$ on the quotient $X$.
\end{sub}

\begin{sub}{\bf Construction of a $G$--linearization of
$\widetilde{M}$:}\label{shmg}\rm

We assume now that $M$ carries a dual $G$--action $M\xrightarrow{\vartheta}
k[G]\otimes_k M$ inducing the dual $H$--action $\vartheta_H$ via the
surjection $k[G]\to k[H]$. Equivalently, see \ref{dact}, we are given a
$G$--linearization
\[
\sigma^\ast \km\xrightarrow{\varphi}p^\ast_2 \km
\]
which restricts to $\varphi_H$ over $H\times Z$. Let $G\times
X\xrightarrow{\widetilde{\sigma}} X$ denote the induced $G$--action on $X$
with diagram
\[
\xymatrix{
G\times W\ar[r]^\sigma\ar[d]_{1\times\pi} & W\ar[d]^\pi\\
G\times X\ar[r]^{\widetilde{\sigma}}& X}
\]
Because $\widetilde{\sigma}$ is flat and the diagram is Cartesian,  
$(1\times \pi)_\ast\sigma^\ast$ and
$\widetilde{\sigma}^\ast \pi_\ast$ (respectively $(1\times \pi)_\ast p_2^\ast$
and $\widetilde{p}_2^\ast \pi_\ast$) are equal functors. Then $\varphi$ yields
a $G$--linearization 
\[
\widetilde{\sigma}^\ast\pi_\ast(\km|W)\xrightarrow{\widetilde{\varphi}}
\widetilde{p}^\ast_2\pi_\ast(\km|W)
\]
together with isomorphisms
\[
g^\ast\pi_\ast(\km|W)\cong\pi_\ast
g^\ast(\km|W)\xrightarrow{\widetilde{\varphi}_g}\pi_\ast(\km|W)
\]
for any $g\in G(k)$.
\end{sub}

\begin{sub}{\bf Claim:} The restriction of $\widetilde{\varphi}$ to
    $\widetilde{\sigma}^\ast\widetilde{\km}_0$ is a $G$--linearization of 
$\widetilde{M}$.
\end{sub}

\begin{sub}\rm
\begin{proof}
By \ref{inv1} it is sufficient to show that $\widetilde{\varphi}_g$ maps
$g^\ast\widetilde{M}$ to $\widetilde{M}$ for any $g\in G(k)$. 
More generally, we show that
$\widetilde{\varphi}_g$ induces an isomorphism
\[
g^\ast\widetilde{\km}_{\chi_g}\xrightarrow{\widetilde{\varphi}_g}\widetilde{\km}_\chi\ ,
\]
where the character $\chi_g$ is defined by $\chi_g(h)=\chi(g^{-1}hg)$. For
that, let $U\subset X$ be open affine, let $\widetilde{U}\subset W$ be its
inverse image in $W$ and consider the diagram
\[
\xymatrix{\widetilde{U}\ar[r]^g\ar[d]^\pi & g\widetilde{U}\ar[d]^\pi\\
U\ar[r]^g & gU}
\]
where $\widetilde{U}$ and $g\widetilde{U}$ are affine and $H$--invariant.
The homomorphism
\[
g^\ast\pi_\ast(\km_W)(U)\xrightarrow{\widetilde{\varphi}_g(U)}\pi_\ast(\km_W)(U)
\]
can be described as 
\[
\ko_Z(\widetilde{U})\otimes_{\ko_Z(g\widetilde{U})} \km(g\widetilde{U})\cong
(g^\ast\km)(\widetilde{U})\xrightarrow{\varphi_g(\widetilde{U})}\km(\widetilde{U})
\]
and this corresponds to a map $\vartheta_g(\widetilde{U})$ as composition of
\[
\km(g\widetilde{U})\to (g^\ast\km)(\widetilde{U})\to\km(\widetilde{U})
\]
with $\varphi_g(\widetilde{U})(s\otimes m)=s\vartheta_g(\widetilde{U})(m)$,
as in \ref{dact}. The cocycle condition corresponds to the commutativity of
the diagrams
\[
\xymatrix{
\km(g_1g_2\widetilde{U})\ar[rr]^{\vartheta_{g_1g_2}(\widetilde{U})}\ar[dr]_{\vartheta_{g_1}(g_2\widetilde{U})}
&& \km(\widetilde{U})\\
& \km(g_2\widetilde{U})\ar[ur]_{\vartheta_{g_2}(\widetilde{U})}}
\]
Using this for a product $hg, h\in H(k), g\in G(k)$, we find that for $m\in
\km(g\widetilde{U})_{\chi_g}$,
\[
\begin{array}{lcl}
\vartheta_h(\widetilde{U})\vartheta_g(\widetilde{U})(m) & = &
\vartheta_{gh}(\widetilde{U})(m)\\
& = & \vartheta_{ghg^{-1}g}(\widetilde{U})(m)\\
& = & \vartheta_g(\widetilde{U})\vartheta_{ghg^{-1}}(g\widetilde{U})(m)\\
& = &
\vartheta_g(\widetilde{U})\chi_g(ghg^{-1})(m)=\chi(h)\vartheta_g(\widetilde{U})(m).
\end{array}
\]
This implies that $\widetilde{\varphi}_g(U)$ maps
$(g^\ast\widetilde{\km}_{\chi_g})(U)$ to $\widetilde{\km}_\chi(U)$ as desired.
\end{proof}
\end{sub}

\begin{sub}{\bf Proposition:}\label{shff} In the above setting,
\begin{enumerate}
\item [(1)] the functor $M\mapsto \widetilde{M}$ is an exact functor from the
  category of $X(H)$--graded $S$--modules to the category of quasi--coherent
  $\ko_X$--modules.
\item [(2)] $\widetilde{M}$ is coherent if $M$ is finitely generated.
\item [(3)] any dual $G$--action $M\xrightarrow{\vartheta} k[G]\otimes M$
over the dual $H$--action on $M$ induces a $G$--linearization 
of $\widetilde{M}$.
\item [(4)] the functors in (1), (2), (3) are essentially surjective.
\end{enumerate}
\end{sub}

\begin{sub}{\bf Remark:}\rm\label{histsh}
Historically, homogeneous coordinates have first been introduced for toric 
varieties, generalizing the standard homogenous coordinates for $\mathbb{P}^n$.
For general toric varieties, these have first been introduced by Cox 
\cite{Cox}.
For simplicial toric varieties it was shown by Cox \cite{Cox}, and by
Musta\c t\v a in general \cite{Mus}, that every quasi-coherent sheaf
on $X$ can be obtained by sheafification of some graded $S$-module. The 
extension
to fine-graded modules and equivariant sheaves was first noted by Batyrev-Cox, 
\cite{BaCo}. Generalizations have been considered in \cite{ahs},
\cite{BerchtoldHausen}, and more recently in \cite{Hausen} for very general 
classes of Cox rings.
\end{sub}

\begin{sub}\rm
\begin{proof} (1), (2), (3) are obvious by construction. To show that
  the functor (1) is surjective up to isomorphisms
  we use the twisted sheaves $\ko_X(\chi_0):=\widetilde{S(\chi_0)}$
  where $S(\chi_0)$ is the shifted $S$--module with
  $S(\chi_0)_\chi=S_{\chi_o+\chi}$ for characters $\chi_0, \chi$ of $H$. Let
  also $\ko_Z(\chi_0)$ denote the sheafification of $S(\chi_0)$ on $Z$. It
  comes with an induced $G$--linearization and thus for any open subset
  $V\subset Z$ with a dual action $\ko_Z(\chi_0)(V)\to
  k[H] \otimes_k \ko_Z(\chi_0)(V)$, inducing an $X(H)$--grading
  $\ko_Z(\chi_0)(V)_\chi$. In fact, for any open affine $U\subset X$, we have
\[
\ko_Z(\chi_0)(\widetilde U)=\oplus_{\chi}\ko_Z(\widetilde{U})_{\chi_0+\chi},
\]
by \ref{shm}. Then, by definition
\[
\ko_X(\chi_0)(U)=\pi_\ast\ko_Z(\chi_0)_0(U)=\ko_Z(\widetilde{U})_{\chi_0},
\]
see notation in \ref{shm}.
\vskip3mm

(i) Let now $\kf$ be a quasi--coherent sheaf on $X$, let
  $\kf(\chi):=\kf\otimes \ko_X(\chi)$ and let
\[
\Gamma_\ast(\kf):=\underset{\chi}{\oplus}\Gamma(X, \kf(\chi))\ .
\]
For any affine open $U\subset X$, also $\widetilde{U}$ is affine and we have
\[
\begin{array}{lcl}
\Gamma(\widetilde{U}, \pi^\ast\kf) & \cong &
\ko_Z(\widetilde{U})\otimes_{\ko_X(U)}\kf(U)\\
 & \cong &
 \underset{\chi}{\oplus}\ko_Z(\widetilde{U})_\chi\otimes_{\ko_X(U)}\kf(U)\\
& \cong &
 \underset{\chi}{\oplus}\ko_X(\chi)(U)\otimes_{\ko_X(U)}\kf(U)\\
& \cong &
 \underset{\chi}{\oplus}\Gamma(U,\kf(\chi))
\end{array}
\]
It follows that $\pi_\ast\pi^\ast\kf\cong
\underset{\chi}{\oplus}\kf(\chi)$. Setting 
$M:=\Gamma(Z, \iota_\ast\pi^\ast\kf)$, where $\iota$ denotes the inclusion of 
$W$ into $Z$, we obtain
\[
\widetilde{M}=(\pi_\ast\pi^\ast\kf)_0\cong \kf\ .
\]
One should note that the grading of $\Gamma(\widetilde{U}, \pi^\ast\kf)$ is
induced by the ''trivial'' dual action arising from $\Gamma(\widetilde{U},
\pi^\ast\kf)\to \Gamma(H\times \widetilde{U}, p_{2H}^\ast\pi^\ast\kf)$.
\vskip3mm

(ii) If $\kf$ is coherent, $\pi^\ast\kf$ is coherent and one can find
  a graded finitely generated submodule $F\subset M$ with
  $\widetilde{F}\cong\kf$.
\vskip3mm

(iii) Finally, suppose that $\kf$ is equipped with a
  $G$--linearization
  $\widetilde{\sigma}^\ast\kf\xrightarrow{\varphi}\widetilde{p}_2^\ast\kf$ on
    $G\times X$ which restricts to the identity on $H\times X$. Pulling this
    back to $G\times W$ and extending to $G\times Z$, we obtain a
    $G$--linearization on $\iota_\ast\pi^\ast\kf$, because
    $\iota_\ast(1\times\pi)^\ast \widetilde{\sigma}^\ast\kf\cong
    \sigma^\ast\iota_\ast\pi^\ast\kf$ (note that $\sigma$ is flat). The
    restriction to $H\times Z$ induces via the dual actions
\[
(\pi^\ast\kf)(\widetilde{U})\to\Gamma(H\times \widetilde{U},
\sigma^\ast_H\pi^\ast\kf)\to k[H]\otimes_k (\pi^\ast\kf)(\widetilde{U})
\]
the grading $\underset{\chi}{\oplus}\Gamma(U, \kf(\chi))$ above because this
is the lift of the trivial dual action on $\kf(U)$. The $G$--linearization of
$\iota_\ast\pi^\ast\kf$ obtained, induces a dual action on $M=\Gamma(Z,
\iota_\ast\pi^\ast  \kf)$. One verifies that this induces the given
$G$--linearization on $\kf\cong\widetilde{M}$.
\end{proof}
\end{sub}

\begin{sub}{\bf Remark:}\rm\label{remdesc} Since the action 
$G\times W\xrightarrow{\sigma} W$ induces also an action $G/H\times X\to X$
it would be desirable to induce a $G/H$--linearization of $\widetilde{M}$ 
by descent. However, in general, such a descent is only possible if
the restriction of the $G$--linearization to $H\times X$ is the 
identity isomorphism, e.g. the tautological subbundle on $\P_n$ does not 
admit a $PGL(n+1)$-linearization. 
In the case of toric varieties however, $\widetilde{M}$ always admits
a $G/H$--linearization as will follow from the more general
\end{sub}

\begin{sub}{\bf Proposition:}\label{existlin}
If $X$ admits a $G$-invariant affine cover, then a dual $G$-action
$M \overset{\vartheta}{\longrightarrow} k[G] \otimes_k M$ induces
a $G/H$-linearization of $\widetilde{M}$.
\end{sub}

\begin{sub}\rm
\begin{proof}
The action of $G / H$ on $X$ is induced by the surjection $G \times X
\rightarrow G / H \times X$. Dually, for every $G$-invariant open affine
subset $U \subseteq X$, we get an injective morphism of rings $k[G / H]
\otimes k[U] \rightarrowtail k[G] \otimes k[U]$, where $k[U]$ denotes
the coordinate ring of $U$. The sheaf $\widetilde{M}$ has a $G$-linearization
by \ref{shmg}. To show that this $G$-linearization induces a
$G / H$-linearization of $\widetilde{M}$,
it suffices to show that on every $U$ the dual action $\Gamma(U, \widetilde{M})
\longrightarrow k[G] \otimes_k \Gamma(U, \widetilde{M})$ restricts to a dual
action $\Gamma(U, \widetilde{M}) \longrightarrow k[G / H] \otimes_k
\Gamma(U, \widetilde{M})$.

For this we can assume without loss of generality that $X$ itself is affine.
The $G$-linearization
corresponds to a dual $G$-action on $\Gamma(X, \widetilde{M}) = M_0$, given by
$M_0 \overset{\vartheta}{\longrightarrow} k[G] \otimes_k M_0$. Explicitly,
$\vartheta$ is given by a mapping $m \mapsto \sum_i\varphi_i \otimes m_i$ for
some $\varphi_i \in k[G]$ and $m_i \in M_0$. Because $m$ is $H$-invariant,
by Lemma \ref{lemdesc},  $\varphi_i \in k[G]^H = k[G/H]$ for all $i$. 
Hence $\vartheta$ restricts
to $M_0 \overset{\vartheta}{\longrightarrow} k[G/H] \otimes_k M_0$. Thus
the $G$-linearization descents to a $G/H$-linearization on $\widetilde{M}$.
\end{proof}
\end{sub}

\begin{sub}{\bf Remark:}\rm\quad
Note that the results of this section do not depend on the characteristic
of $k$. However, if $H$ has torsion which is not coprime to 
$\operatorname{char} k$, one should
pass to the group scheme $\Spec(k[X(H)])$ rather than the group $H$ itself.
\end{sub}

\begin{sub}{\bf Remark:}\rm\quad
In general it is not true that $\Gamma_* \mathcal{O}_X = S$. However, this 
holds
if we assume that $Z$ is normal and the codimension of $W$ in $Z$ is at least 
two.
\end{sub}

\section{Homogeneous coordinate rings and primary decomposition}\label{quotientprimdec}


In this section we assume the same setting as in section \ref{sqp}. That is, 
$H$ is a diagonalizable normal subgroup of an algebraic group $G$ which acts 
on $Z$,
such that $W \subseteq Z$ is $G$-invariant. Moreover, we assume, as in 
Proposition
\ref{existlin}, that $X$ (and thus $W$) admits a $G$-invariant affine cover.
We denote $V$ the complement of $W$ in $Z$, which itself is $G$-invariant and 
thus
$H$-invariant.

\begin{sub}\rm{\bf $G$-equivariant primary decomposition and graded modules}\label{eprimdec}\\
Let $\kf \subseteq \ke$ be $G$-equivariant sheaves on $Z$ and
$F \subseteq E$ their corresponding $S$-modules. By Theorem \ref{ncepd} we have
a $G$-equivariant primary decomposition $\kf = \underset{i \in I}{\cap} \kq_i$
in $\ke$ with $G$-equivariant sheaves $\kq_i$ which are $Y_i$-primary with
associated $G$-components $Y_i$ of $Z$. The $\kq_i$ also correspond to $S$-modules
$Q_i \subseteq E$. By \ref{dact}, all these modules are equipped with a dual $G$-action.
In particular, all modules are als endowed with a dual $H$-action and thus 
with a
$X(H)$-grading by Corollary \ref{dactdiagcor}. Therefore we get a decomposition
of $X(H)$-graded modules $F = \underset{i \in I}{\cap} Q_i$.
This decomposition is also compatible with the dual $G$-actions. The $G$-components
$Y_i$ are $H$-invariant and thus correspond to homogeneous ideals $\mathfrak{p}_i$
of $S$ with respect to the $X(H)$-grading.

Note that in the sequel we will call a decomposition $F = \underset{i \in I}{\cap} Q_i$
a {\em $G$-equivariant} primary decomposition if the corresponding decomposition
$\kf = \underset{i \in I}{\cap} \kq_i$ is a $G$-equivariant primary decomposition.
\end{sub}

\begin{sub}\rm{\bf Graded primary decompositions}\\
Bourbaki's Commutative Algebra \cite{NBC} treats graded primary decomposition
for the case where the grading is given by a torsion free abelian group.
This has been generalized
in \cite{pk} to the case where the grading may have torsion. For this, one introduces
analogously to \ref{assgc} a slightly more general notion of associated ideals
and primary modules.
\end{sub}

\begin{sub}{\bf Definition:}\label{gprimdecdef}
Let $A$ a finitely generated abelian group, $R$ a noetherian $A$-graded ring,
and $E$ a finitely generated $A$-graded $R$-module.
\begin{enumerate}[(i)]
\item A graded ideal $\mathfrak{p} \subset R$ is $A$-prime if for every two
{\em homogeneous}
elements $r,s \in R$ holds that $rs \in \mathfrak{p}$ implies that $r \in 
\mathfrak{p}$
or $s \in \mathfrak{p}$.
\item An ideal $\mathfrak{p}$ of $R$ is $A$-associated iff it is $A$-prime and
$\mathfrak{p} = \operatorname{ann}(x)$ for some element $x \in E$.
\item $\operatorname{Ass}^A(E)$ denotes the set of all $A$-associated ideals 
of $E$.
\item An $A$-graded submodule $F$ of $E$ is said to be $A$-primary if the 
quotient
module $\operatorname{Ass}^A(E/F) = \{\mathfrak{p}\}$ for
some $A$-prime ideal $\mathfrak{p}$.
\item For $F$ an $A$-graded submodule of $E$, we call an expression $F =
\underset{i \in I}{\cap} Q_i$ an $A$-primary decomposition of $F$ in $E$
iff the $Q_i$ are $A$-primary submodules of $E$ with 
$\operatorname{Ass}^A(E/Q_i)
= \{\mathfrak{p}_i\}$ and the $\mathfrak{p}_i$ are $A$-associated to $E/F$.
\item An $A$-primary decomposition $F = \underset{i \in I}{\cap} Q_i$ of $F$
in $E$ is called reduced if all the $\mathfrak{p}_i$ are
distinct and there exists no $i \in I$ such that $\underset{j \neq i}{\cap} Q_j
\subset Q_i$.
\end{enumerate}
\end{sub}

With respect to these notions, the following results have been proved in 
\cite{pk}:

\begin{sub}{\bf Proposition:}\label{grdec}
({\rm\cite[Theorem 1.2]{pk}})
Let $A$ be a finitely generated abelian group, $R$ a noetherian $A$-graded ring,
$E$ a finitely generated $A$-graded $R$-module, and $F$ an $A$-graded submodule
of $E$. Let $F = \underset{i \in I}{\cap} Q_i$ be a primary decomposition of 
$F$
in $E$. Then:
\begin{enumerate}[(i)]
\item Let $Q_i'$ be the largest $A$-graded submodule of $E$ contained in $Q_i$.
Then $Q_i'$ is $A$-primary for every $i \in I$ and $F = \underset{i \in 
I}{\cap} Q_i'$.
\item There exists a subset $J$ of $I$ such that $F = \underset{i \in J}{\cap} 
Q_i'$
is a reduced $A$-primary decomposition.
\item If some $Q_i$ corresponds to an $A$-prime ideal $\mathfrak{p}_i$ which 
is a
minimal element of $\operatorname{Ass}^A(E/F)$, then $Q_i$ is $A$-graded.
\end{enumerate}
\end{sub}

\comment{
The following shows that $A$-primary decomposition preserves all information 
about the irreducible branches:

\begin{sub}{\bf Proposition:}\label{grdecbranches}({\rm\cite[Corollary 4.5]{pk}})
Let $A$, $R$, $E$, $F$ be as before, except that $R$ is a $k$-algebra. Let
$\operatorname{Ass}^A(E/F) = \{\mathfrak{p}_i\}_{i \in I}$ and consider some
$A$-primary decomposition $F = \underset{i \in I}{\cap} Q_i$, where the
$Q_i$ are $A$-primary submodules of $E$ with respect to $\mathfrak{p}_i$.
For every $i \in I$ denote $\operatorname{Ass} R/\mathfrak{p}_i =: 
\{\mathfrak{p}_{ij}
\}_{j \in J_i}$. Then:
\begin{enumerate}[(i)]
\item $\operatorname{Ass}(E/F)= \underset{i \in I}{\cup} \operatorname{Ass}
A/\mathfrak{p}_i$.
\item For every $Q_i$ we have $\operatorname{Ass}(E/Q_i) = \underset{j \in 
J_i}{\cup}
\operatorname{Ass} A/\mathfrak{p}_{ij}$ and $Q_i$ has a non-graded primary
decomposition $Q_i = \underset{j \in J_i}{\cap} Q_{ij}$, where 
$\operatorname{Ass}
E/Q_{ij}$ = $\{\mathfrak{p}_{ij}\}$.
\end{enumerate}
\end{sub}
}

\begin{sub}\rm{\bf $G$-equivariant decomposition on $X$}\\
It follows straightforwardly that the $G$-equivariant primary
decomposition of \ref{eprimdec} induces a $X(H)$-primary
decomposition in the sense of Definition \ref{gprimdecdef}.
Denote $V$ the complement of $W$ in $Z$ and $B$ the radical
ideal in $S$ whose zero set is $V$. This ideal is
$X(H)$-homogeneous and we call it the {\em irrelevant ideal}.
If $B \subset\mathfrak{p}_i$ for some $i \in I$, then the support of 
the module $Q_i/F$ is contained in $V$, and thus its sheafification
necessarily vanishes.
Denote $\kf[V]$ the relative gap sheaf and
$F[V] = \Gamma(Z, \kf[V])$.
By Propositions \ref{gapdec} and \ref{grdec}
we get the corresponding primary decompositions
\begin{equation*}
\kf[V] = \underset{Y_i \not \subset V}{\cap} \kq_i \ \text{ and }\ 
F[V] = \underset{B \not \subset \mathfrak{p}_i}{\cap} Q_i,
\end{equation*}
where the latter is also a $X(H)$-primary decomposition.
 As the sheafification
is compatible with taking intersections, it follows that
$$\widetilde{F} \cong \widetilde{F[V]}.$$
In the special case where $F = 0$, the relative gap sheaf $0[V]$
coincides with the local cohomology module $H^0_B E$, and consequently
$\widetilde{0[V]} = 0.$
\end{sub}
\vskip3mm

\begin{sub}{\bf Proposition:}\label{tordec1}
Let $F \subseteq E$ be finitely generated $X(H)$-graded $S$-modules. Then 
there exists a decomposition of sheaves
$$
\widetilde{F} = \underset{\substack{i \in I\\ B \not\subset \mathfrak{p}_i}}{\cap}
\widetilde{Q}_i\subset \widetilde{E}.
$$
on $X$. If the $E$, $F$ admit a dual $G$-action, then under the 
assumption on the action of $G$, this decomposition
is also $G/H$-equivariant.
\end{sub}

\begin{proof}
The assertions follow by choosing a $G$-equivariant primary decomposition for
$\kf$ or, equivalently, a $X(H)$-graded primary decomposition with gradings
induced by dual actions of $G$. The sheaves 
$\widetilde{F}$ and $\widetilde{Q}_i$
are $G/H$-equivariant by Proposition \ref{existlin}.
\end{proof}

\begin{sub}\rm
So far, it is not clear whether the decomposition of Proposition \ref{tordec1} is
a global $G/H$-equivariant primary decomposition in the sense of section \ref{sectepd}.
The sheaves $\widetilde{Q}_i$ may not be $G/H$-primary sheaves on $X$. To
establish this property we have to require that $H$ acts freely on $W$.
Then for every point $w \in W$ the image of the morphism
$H \rightarrow H . z$, given by $h \mapsto h . z$, is isomorphic to $H$.
In particular, $X = W // H$ is a geometric quotient. 
Note that in general for $W // H$ being a geometric quotient, the action
of $H$ on $W$ is allowed to have finite stabilizers. In the presence of
nontrivial stabilizers it is not true in general that the decomposition
of Proposition \ref{tordec1} induces a primary decomposition of sheaves on
$X$, as we will see in examples \ref{singexample1} and
\ref{singexample2} below.
\end{sub}


\begin{sub}{\bf Lemma:}\label{lemdeg}
Assume that $H$ acts freely on $W$.
Then every  $x \in X$ has an affine neighbourhood $U_x$
such that $\pi^{-1}(U_x)$ admits invertible homogeneous
functions in any degree and we have $\pi^{-1}(U_x) = Z_g =
\{x \in Z \mid g(x) \neq 0\}$ for some homogeneous $g \in S$.
\end{sub}

\begin{proof}
Fix a point $z \in \pi^{-1}(x)$ and choose some homogeneous $f \in B$ such that
$z \in Z_f$, where $Z_f = \{x \in Z \mid f(x) \neq 0\}$. We denote $S_f$ the
coordinate ring of $Z_f$, which is the localization of $S$ by $f$. By the fact
that the quotient $\pi : W \rightarrow X$ is geometric, we obtain on the one
hand that
the orbits $H . z$ are closed in $Z_f$ and therefore any homogeneous function
on $H . z$ extends to a homogeneous function on $S_f$. Because $H$ acts freely on
$W$, the orbit $H . z$ is isomorphic to $H$. Therefore we
can identify the coordinate ring of $H . z$ in a natural way with the group
ring $k[X(H)]$, and, moreover, we can lift every character
$\chi \in k[X(H)]$ to some homogeneous element in $(S_f)_\chi$. In particular, for
any system of generators $\chi_1, \dots, \chi_r$ of $X(H)$, we can choose $g_i \in
(S_f)_{\chi_i}$ with $g_i(z) \neq 0$. So we obtain an affine subset $W_x :=
Z_f \bigcap_i Z_{g_i} = Z_g$, where $g = f^t g_1, \cdots g_r \in S$ for $n$
big enough, and $U_x := \pi(W_x)$, as wanted.
\end{proof}

\comment{
\begin{sub}{\bf Lemma:}\label{lemdeg}
Let $X$ be a smooth toric variety and $S$ its homogeneous coordinate ring. 
Then,
for every $\sigma \in \Delta$, we consider the localization $S_{x(\sigma)}$.
In this ring, every homogeneous element $f \in (S_{x(\sigma)})_\alpha$ can be
factorized as $f = f' \cdot u$ where $f'$ is of degree zero and $u$ is a unit.
\end{sub}

\begin{proof}
Recall that $x(\sigma) = \prod_{\rho \in \rays \setminus \sigma(1)} x_\rho$,
which implies that every $x_\rho$, $\rho \in \rays \setminus \sigma(1)$ becomes
a unit in $S_{x(\sigma)}$. First, we choose any representative 
$a = (a_\rho \mid \rho \in \rays)
\in \mathbb{Z}^\rays$ of $\alpha \in X(H)$. 
Without loss of generality, we may assume
that $a_\rho > 0$ for all $\rho \in \sigma(1)$. Else, we can choose some $m \in
\sigma_M$ such that $\langle m, n(\rho) \rangle > 0$ for every 
$\rho \in \sigma(1)$;
adding a sufficient high multiple of the image of $m$ in $\mathbb{Z}^\rays$ to
$a$ then gives some $a'$ such that $a'_\rho > 0$ for all
$\rho \in \sigma(1)$ which also is a representative of $\alpha$. Because $X$ is
smooth, the primitive vectors $n(\rho)$, for $\rho \in \sigma(1)$ are part of
a $\mathbb{Z}$-basis of $F$. Therefore the system of equations 
$\langle m, n(\rho)\rangle = a_\rho$ 
has at least one integral solution $m_0$. So, denoting $\hat{m}_0
:= (\langle m_0, n(\rho) \rangle \mid \rho \in \rays)$, we can write
$a = \hat{m}_0 + b$. The monomial $x^b$ is a
unit in $S_{x(\sigma)}$ and has degree $\alpha$. We set
$u := x^b$ and $f' := f \cdot x^{-b}$.
\end{proof}
}

We use this lemma to prove:

\begin{sub}{\bf Theorem:}\label{tordec2}
Assume that $H$ acts freely on $W$.
Let $F\subset E$ be finitely generated $G$-equivariant $S$-modules and let 
$F = \underset{i \in I}{\cap} Q_i$
be an $G$-equivariant primary decomposition of $F$ in $E$. Then the decomposition
$$
\widetilde{F} = \underset{\substack{i \in I\\ B \not\subset \mathfrak{p}_i}}{\cap}
\widetilde{Q}_i\subset \widetilde{E}
$$
is a global $G/H$-equivariant primary decomposition on $X$.
\end{sub}

\begin{proof}
The $G / H$-equivariance follows from Proposition \ref{existlin}. So it
remains to show that $\widetilde{E} / \widetilde{Q}_i$ is $Y$-primary
for some associated $G/H$-component $Y$ on $X$. Let $x$ be any point in $X$.
By Lemma \ref{lemdeg} there exists an affine neighbourhood $U$ of $x$
such that $\pi^{-1}(U) =: \widetilde{U}$ is affine and whose coordinate
ring, which we denote $S_U$, is $X(H)$-graded and admits invertible
homogeneous functions in any degree. This implies that every $X(H)$-graded
$S_U$-module is generated in degree $0$ and thus there is a
one-to-one correspondence between $(S_U)_0$-modules and $X(H)$-graded
$S_U$-modules which is given by taking degree zero. In particular, there
is a one-to-one correspondence of homogeneous ideals in $S_U$ and ideals
in $(S_U)_0$. We can interpret this geometrically as the correspondence
between subschemes of $U$ and $H$-invariant subschemes of $\tilde{U}$,
where the associated scheme $V(\mathfrak{a}) \subseteq \widetilde{U}$
of a homogeneous ideal $\mathfrak{a}$ corresponds to $V(\mathfrak{a}_0)
\subseteq U$. The geometric
quotient $U = \widetilde{U} // G$ restricts to a geometric quotient
$\pi|_{V(\mathfrak{a})}: V(\mathfrak{a}) \longrightarrow
V(\mathfrak{a}_0) = V(\mathfrak{a}) // G$ with constant fiber dimension
$\kappa = \dim H$.

Now let $M$ be any $X(H)$-graded $S_U$-module. Any homogeneous unit
$u \in S_h$ for some $h \in H$, yields a bijection $M_0
\overset{\cdot u}{\longrightarrow} M_h$. It is straightforward to see
that this bijection respects supports, i.e. for any $d \geq 0$ we have
an induced bijection $\Gamma(\widetilde{U}, \kh_d^0\mathcal{M})_0
\overset{\cdot u}{\longrightarrow} \Gamma(\widetilde{U}, \kh_d^0\mathcal{M})_h
\subset M$, where $\mathcal{M}$ is the sheaf on $\widetilde{U}$
associated to $M$. This implies that there is a bijection between
the associated $H$-components of $\mathcal{M}$ and the associated
components of $\widetilde{M}$. Now assume that $\mathcal{M}$
additionally is $G$-equivariant; then $Y \subseteq \Supp
\kh^0_d \mathcal{M}$ is an associated $G$-component of $\mathcal{M}$
if and only if $Y // G \subseteq \Supp \kh^0_{d - \kappa} \widetilde{M}$
is an associated $G / H$-component of $\widetilde{M}$.

Hence, $\widetilde{F} =
\underset{\substack{i \in I\\ B \not\subset \mathfrak{p}_i}}{\cap}
\widetilde{Q}_i\subset \widetilde{E}$ is a $G/H$-equivariant
primary decomposition.
\end{proof}

\section{Examples for equivariant primary decompositions on toric varieties}\label{toricsh}

In this section we consider examples of equivariant primary decompositions over
toric varieties as a special case of examples for the setting of sections 
\ref{sqp} and \ref{quotientprimdec}. A toric variety over $k$ is a normal
irreducible variety $X$ acted on by the torus $T=(k^{\ast})^n$, where
$n=\dim(X)$, in such a way that there is an open dense orbit which is
isomorphic to $T$ and such that the action on this orbit is the multiplication
of $T$. We briefly recall some basic notation, definitions and facts for toric 
varieties and sheaves on them, which are used in the description of primary
decompositions. As basic references for the theory of toric varieties we refer 
to \cite{Od} and \cite{Fu}.
\vskip3mm 

\begin{sub}{\bf Toric varieties and divisors}\label{tordiv}\\ \rm
Let $M=X(T)\cong \Z^n$ be the free character group of the torus and $N$
denote its dual. The toric variety is characterized by a finite fan $\Delta$
of strongly  convex rational polyhedral cones $\sigma\subset N_\R =
N\otimes_\Z \R\cong \R^n$. For simplicity, we will assume in the sequel
that the fan is not contained in a proper subvector space of $N_\R$.
Any cone $\sigma\in\Delta$ determines an open affine $T$-invariant
subset $U_\sigma\subset X$ and a closed irreducible $T$-invariant 
subvariety $V(\sigma)\subset X$, the closure of the minimal orbit
in $U_\sigma$, such that $\dim V(\sigma)= n - \dim\, \sigma.$ 
We denote $\Delta(1)$ and $\sigma(1)$, respectively, the set of
$1$-dimensional cones $\rho\in\Delta$ and $\rho\in\sigma$,
respectively. Each one-dimensional cone $\rho$
is spanned by a primitive lattice vector $n(\rho)\in\rho\cap N$
and determines an irreducible T-invariant divisor $D_\rho = V(\rho).$

For any $m\in M$ let $\chi (m)$ denote the character function with the
rule $\chi (m+m')= \chi (m)\chi (m').$ Then $\chi (m)$ is a rational
function on $X.$  A proof of the following statements is given in \cite{Fu}.

(1) The group of $T$-invariant Weil divisors of $X$ is freely
generated by the divisors $D_\rho.$ We will denote this group
$\Z^\rays$ and we write $D(a)=\Sigma a_\rho D_\rho$ for the divisor 
corresponding to $a\in\Z^\rays$. 

(2) For any $m\in M$ the divisor of $\chi (m)$ is given by
$$\dv\chi(m)= \underset{\rho\in\Delta(1)}\Sigma 
\langle m, n(\rho) \rangle D_\rho .$$

(3) The divisor class group is generated by the $T$-invariant Weil divisors.
More precisely, there is a short exact sequence 
\begin{equation*}
0 \lra M \overset{div}\lra \mathbb{Z}^\rays \overset{cl}\lra A_{n - 1}(X) \lra 
0.
\end{equation*}
Here, $cl$ denotes the class map to the divisor class group of $X$. From now,
we will abbreviate $A_{n - 1}(X)$ by $A$.
\end{sub}
\vskip3mm

\begin{sub}{\bf Homogeneous coordinates for toric varieties}\label{hcori}\\ \rm
To specialize the setting of section \ref{quotientprimdec} to the case of toric
varieties, we consider the vector space $k^\rays = \mathbb{Z}^\rays
\otimes_\mathbb{Z} k$ which itself is an affine toric variety. Its associated
cone can be identified in a natural way with the positive
orthant of the vector space $\R^\rays$. Its primitive vectors then are given
by the standard basis $\{e_\rho\}_{\rho \in \rays}$ of $\R^\rays$. We consider
a subfan $\hat{\Delta}$ of this cone which is generated by cones 
$\hat{\sigma}$,
where $\hat{\sigma} \in \hat{\Delta}$ whenever there exists some $\sigma'
\in \Delta$ such that the set $\{\rho \mid e_\rho \in
\hat{\sigma}(1)\}$ is contained in $\sigma'(1)$.

The fan $\hat{\Delta}$ defines an open dense toric subvariety $\hat{X}$ of
$k^\rays$ together with a natural toric morphism $\pi: \hat{X} \longrightarrow 
X$.
The latter is given by the map of fans which is induced by mapping $e_\rho$
to $n(\rho)$ for every $\rho \in \rays$. This morphism is equivariant with
respect to the actions of the tori $T = \Hom(M, k^*)$ and $\hat{T} = 
\Hom(\Z^\rays, k^*)$.
We obtain a short exact sequence of groups
\begin{equation*}
1\lra H \lra \hat{T} \lra T \lra 1,
\end{equation*}
where $H = \Hom(A, k^*)$ is a diagonalizable subgroup of the torus $\hat{T}$
with $X(H) = A$. In fact,
this short exact sequence is nothing but the divisor class sequence above 
dualized
by $\Hom(\, . \, , k^*)$. It was shown in \cite{Cox} that in fact $X$ is a good
quotient of $\hat{X}$ by $H$ which is geometric if and only if $\Delta$ is a
simplicial fan. By the short exact sequence and the fact that $\hat{X}$ has
an affine $\hat{T}$-invariant cover, the conditions of \ref{existlin} are
satisfied. In particular, the homogeneous coordinate ring $S$ in this situation
is the polynomial ring
\begin{equation*}
S := k[x_\rho \mid \rho \in \rays],
\end{equation*}
and inherits two natural gradings. The first is the $\Z^\rays$- or 
fine-grading
coming from the dual action of $\hat{T}$ on $S$, where $\deg_{\Z^\rays}(x^a)
= a$ for any $a \in \N^\rays$. The second is the $A$-grading given by
$\deg_A(x^a) = cl(a) \in A$ for every $a \in \N^\rays$. Corresponding to
these gradings we write
$$S=\oplus_{a\in\Z^\rays}S_a =\oplus_{\alpha\in A}S_\alpha\quad\text{with}\quad
S_\alpha=\oplus_{cl(a)=\alpha}S_a.$$

The complement of $\hat{X}$ in $k^\rays$ is described by the irrelevant ideal
$B\subset S$, which is generated by the monomials $x(\sigma),\, 
\sigma\in\Delta$,
where $x(\sigma) := \underset{\rho\not\in \sigma(1)}\Pi x_\rho.$

For any of these monomials one has the $A$-graded localization 
$S_{x(\sigma)}$ and it turns out 
that the subring of elements of $A$-degree zero, 
$$S^\sigma:= S_{(x(\sigma))}= (S_{x(\sigma)})_0$$  
is isomorphic to the coordinate ring 
$k[U_\sigma]$, see \cite{Cox}.
Under this isomorphism the monomial in $S^\sigma$ corresponding to the 
character $\chi(m)$ is 
$$x^m:=\underset{\rho \in \rays}\Pi x_\rho^{\langle m, n(\rho) \rangle}\, .$$
\end{sub}

\begin{sub}{\bf The associated sheaves:}\label{toricsh2}\quad \rm
The sheafification of Section \ref{sqp} specializes in the above 
toric situation as follows, as originally described in \cite{Cox}. Given 
an $A$-graded $S$-module $F$, there are the localized 
$S_{(x(\sigma))}=k[U_\sigma]$-modules $F^\sigma= F_{(x(\sigma))}$ so that the sheaf  
$\widetilde{F}$ is determined by 
$\widetilde{F}|U_\sigma= \widetilde{F}_{(x(\sigma))}.$
\vskip2mm

When $F$ is a fine-graded $S$-module, the sheaf $\widetilde{F}$ comes 
with an induced $T$-linearization 
$\theta^\ast \widetilde{F}\xrightarrow{\Phi} p_2^\ast\widetilde{F},$
where $T\times X\xrightarrow{\theta} X$ denotes the torus action,
cf. \ref{existlin}. In this case the localized module 
$F^\sigma = F_{(x(\sigma))}$
is an $M$--graded $k[U_\sigma]$--module with components
\[
F^\sigma_m=\{f/x(\sigma)^k\ |\ f_a\in F_a,\, a\in \Z^{\rays},\, 
a-ka(\sigma)=m\in M\}.
\]

This $M$-grading defines naturally a dual $T$-action on 
$F^\sigma,$ see \ref{dactdiagcor}, and thus a $T$-linearization 
\[
\theta^\ast \widetilde{F}^\sigma\xrightarrow{\Phi_\sigma} p_2^\ast
\widetilde{F}^\sigma
\]
over $T\times U_\sigma$. 
The isomorphisms $\Phi_\sigma$ satisfy the cocycle condition and coincide on
intersections, thus defining the global $T$-linearization of $\widetilde{F}$.
\end{sub}

\begin{sub}{\bf The divisorial sheaves:}\label{divisorial}\\ \rm
For any $\alpha\in A$ there is the twisted module $S(\alpha)$ defined by 
$S(\alpha)_\beta=S_{\alpha+\beta}.$  Its associated sheaf 
$\widetilde{S(\alpha)}$ is denoted by $\ko_X(\alpha)$.
\vskip2mm

Similarly, for any $a \in \Z^\rays$ there is the fine-graded twisted 
module $S(a).$
Its associated sheaf  $\widetilde{S(a)}$ on $X$ is canonically isomorphic to
the divisorial sheaf $\ko_X(D(a)).$ It turns out that 
the canonical $T$-linearization of $\widetilde{S(a)}$ corresponds to that of 
$\ko_X(D(a))$ as defined in \ref{equidiv}. 

Furthermore, for effective divisors $D(a)$ the spaces 
$\Gamma(X,\ko_X (D(a))^T$ of $T$-invariant sections are 1-dimensional
with $$\Gamma(X,\ko(D(a)))^T=S_a=k\cdot x^a$$
as in the case  of $\P_n$ in \ref{torlb}.

The sheaves $\ko(D(a))$ and $\ko(\alpha)$ are isomorphic whenever
$cl(a)=\alpha.$ 

It has been proved in \cite{Cox} that $\Gamma(X,\ko_X(\alpha))\cong S_\alpha$.
If we choose a $T$-invariant representative $D(a), a \in \Z^\rays$, then we 
obtain a natural isotypical decomposition
\begin{equation*}
\Gamma(X,\ko_X(D(a)) \cong \bigoplus_{m \in M} 
\Gamma(X,\ko_X(D(a))_m
= \bigoplus_{a \in \N^\rays, cl(a) = \alpha} S_a.
\end{equation*}

\end{sub}

\begin{sub}{\bf Equivariant matrices:}\label{matcd}\\ \rm
Choosing the canonical T-linearizations of the sheaves $\ko_X(D(a)),$
we obtain for any two $a, b \in \Z^\rays$ a canonical isomorphism
\[
\Hom(X, \ko(D(a)), \ko(D(b)))^T\cong\Gamma(X, \ko(D(b)-D(a)))^T
\]
between the space of $T$--invariant homomorphisms and the 
space of $T$--invariant sections. Each of them is either $0$-- or 
$1$--dimensional.

So, if $a\leq b$, there is, up to scalar multiplication, a unique 
$T$-equivariant
homomorphism which can naturally be written as $\ko(D(a))\xrightarrow{x^{b-a}}
\ko(D(b))$. Now, a general homomorphism
\[
\underset{i}{\oplus} \ko (\alpha_i) \overset{F}{\lra} \underset{j}{\oplus} 
\ko(\beta_j)
\]
with $\alpha_i, \beta_j \in A$ is given by a matrix $F$ whose entries are
$A$-homogeneous polynomials in $S$. Such a homomorphism is $T$-linearizable
if there exist $a_i, b_j \in \Z^\rays$ with $cl(a_i) = \alpha_i$, $cl(b_j) = 
\beta_j$ and
associated toric divisors $D(a_i)$, $D(b_i)$ such that
$F$ is induced by a $T$--equivariant homomorphism
\[
\underset{i}{\oplus} \ko\big(D(a_i)\big) \overset{G}{\lra} 
\underset{j}{\oplus} \ko
\big(D(b_j)\big)\,.
\]
In that case, $G$ consists of monomials and $\Coker(F)$ admits 
a $T$--linearization. Moreover, in such a matrix all minors are
monomials, too, as can immediately be verified.   
Conversely, this property characterizes the torus invariance:
\end{sub}
\vskip3mm

\begin{sub}{\bf Proposition:}\label{invh}\quad Let $F = (f_{ij})$ be a 
matrix of monomials\\  
$f_{ij} = \lambda_{ij} x^{a_{ij}}\in\Gamma(X, \ko(\beta_j - \alpha_i)).$
Then the following conditions are equivalent:
\begin{enumerate}
\item [(1)] $F$ is $T$--linearizable
\item [(2)] All minors of $F$ are monomials.
\item [(3)] There are $a_i, b_j \in \Z^{\Delta(1)}$ such that 
$a_{ij} = b_j - a_i$, whenever $f_{ij} \not= 0$.
\end{enumerate} 
\end{sub}

This can be proved by elementary matrix operations using induction on the 
size of the matrices. For a proof in the category of graded modules 
see \cite{CD}.

\begin{sub}{\bf Remarks:}\label{isot} \rm
(1) Homogeneous coordinate rings are in general far from unique. For the toric 
case, other examples
for such presentations have been constructed by Kajiwara \cite{Kajiwara} (see 
also
\cite{P5}). However, the standard Cox construction seems to be the most simple
and best suited for our purposes. For general properties of homogeneous 
coordinate rings
in the toric case see \cite{ahs}.\\
(3) Note that for any equivariant and coherent sheaf $\kf$, the  isotypical
component $\Gamma(U_\sigma, \kf)_m$ of the $T$--module $\Gamma(U_\sigma, \kf)$
is finite dimensional (see \cite{P1}, \S 5.3).
\end{sub}

\begin{sub}{\bf The sheaf of Zariski differential forms}
\label{zariskisheaves}\\ \rm
Denoting $U \overset{i}{\hookrightarrow} X$ the imbedding of the open subset of
regular points of $X$, the sheaves $\widetilde{\Omega}^p_X$ of Zariski 
differential
forms are defined as the reflexive extension $i_\ast \Omega^p_U$ of the usual 
sheaves
of differentials on $U$.
$\widetilde{\Omega}^p_X$ is the kernel of the first
homomorphism of the Ishida--complex, see \cite {Od}, Ch. 3, i.e. there is an
exact sequence
\[
0\to \widetilde{\Omega}^p_X \to \Lambda^p M\otimes_\Z\ko_X\to
\underset{\rho\in\Delta(1)}{\oplus} \Lambda^{p-1}(\rho^\perp\cap
M)\otimes_\Z\ko_{D_\rho}
\]
for any $p$. For $p=1$ this sequence reduces to
\[
0\to \widetilde{\Omega}^1_X\to
M\otimes_\Z\ko_X\xrightarrow{\beta}\underset{\rho\in
  \Delta(1)}{\oplus}\ko_{D_\rho}\to\kc\to 0
\]
where $\kc$ denotes the cokernel of the homomorphism $\beta$ which is
naturally defined as the composition $M\otimes \ko_X\to \Z^{\Delta(1)}\otimes
\ko_X\to \underset{\rho}{\oplus}\ko_{D_\rho}$ of the inclusion
$M\xrightarrow{\beta}\Z^{\Delta(1)}$ and the direct sum of the restrictions
$\ko_X\to \ko_{D_\rho}$.

The embedding of the sheaf $\widetilde{\Omega}^1_X$ into
$M\otimes_\Z \ko_X$ provides an example of an equivariant primary decomposition
on a toric variety. For this, we first have to take a closer look at the toric
Euler sequence. Following \cite{BaCo}, we use the divisor class sequence
\ref{tordiv} and
obtain the commutative exact diagram
\[
\UseComputerModernTips
\xymatrix{
 & 0\ar[d] & 0\ar[d] & & &\\
0\ar[r] & \widetilde{\Omega}^1\ar[d]\ar[r] &
M\otimes_\Z\ko_X\ar[r]\ar[d] &
\oplus_\rho\ko_{D_\rho}\ar[r]\ar@{=}[d] &\kc\ar[r]  & 0\\
0\ar[r] &  \oplus_\rho \ko_X(-D_\rho)\ar[r]\ar[d] &
\Z^{\Delta(1)}\otimes_\Z \ko_X\ar[r]\ar[d] &\oplus_\rho\ko_{D_\rho}\ar[r] & 0\\
 & A\otimes_\Z\ko_X\ar[d]\ar@{=}[r] & A\otimes_\Z\ko_X \ar[d]& & \\
 & \kc\ar[d] & 0 & & &\\
 & 0 & & & &}
\]
The vertical left column generalizes the well-known Euler sequence on a
projective space and has been used in \cite{BaCo} for simplicial toric 
varieties.
All the sheaves in this
diagram are $T$--equivariant with $T$--invariant homomorphisms because it is
the sheafification of the fine--graded diagram
\[
\UseComputerModernTips
\xymatrix{
 & 0\ar[d] & 0\ar[d] & & &\\
0\ar[r] & F\ar[r]\ar[d] & S^n\ar[r]^(0.3){\beta}\ar[d]^Q & 
\oplus_\rho S/S(-e_\rho)\ar[r]\ar@{=}[d] & C\ar[r] & 0\\
0\ar[r] & \oplus_\rho S(-e_\rho)\ar[r]\ar[d]_\alpha & 
S^r\ar[r]\ar[d]^P &\oplus_\rho S/S(-e_\rho)\ar[r] & 0\\
 & S^{r-n}\ar@{=}[r]\ar[d] & S^{r-n} \ar[d]& & &\\
 & C\ar[d] & 0 & & &\\
& 0 & & & &
}
\]
For that we have chosen a $\Z$--basis of $M$ and of $A$/torsion such that
the matrices $Q$ and $P$ describe the divisor class sequence and have integers
as entries. Finally, $\beta$ and $\alpha$ are the compositions
in the diagram and respect the fine--grading. Note that $\alpha$ is the
product
\[
P\circ
\left(\begin{array}{ccc}
x_1 & & \\
& \ddots & \\
 & & x_r
\end{array}\right)
\]
so that all its minors are monomials again.
\end{sub}

It had been shown in \cite{Od}, ch.3, and \cite{Da} that $\alpha$ or $\beta$
is surjective over $U_\sigma$ if $\sigma$ is simplicial, see also
\cite{BaCo}. The converse also holds:

\begin{sub}{\bf Lemma:}\label{sipl}\; 
$\kc|U_\sigma=0$ if and only if $\sigma$ is simplicial.

\begin{proof}
Set $E:=S^n$,\quad $Q_\rho:=S/S(-e_\rho),$ and 
 $F^\sigma : = \Gamma(U_\sigma, \widetilde{F})$ 
for any $A$-graded $S$-module $F$. When $F$ is fine-graded, $F^\sigma_m$ 
denotes the component of the $M$-graded module $F^\sigma$ of degree $m$. 
Moreover, we conveniently identify
$M$ with its image $div(M)$ in $\Z^\rays$, see \ref{tordiv}, 
and use the notation $\sigma_M:= \sigma^\vee\cap M,$ where $\sigma^\vee$
is the dual of the cone $\sigma.$ 

The sheaf $\kc$ vanishes on $U_\sigma$ if and only if the induced maps
$$\beta_m^\sigma\, : E_m^\sigma\to (Q_\rho)_m^\sigma $$
are surjective for any $m\in M$. Now a direct verification shows that
\[
S^\sigma_m\cong
\left\{
\begin{array}{ll}
k\cdot x^m & \text{if}\quad m\in \sigma_M\\
0 & \text{otherwise}
\end{array}\right .
\quad \text{ and }\quad  S(-e_\rho)^\sigma_m\cong
\left\{
\begin{array}{ll}
k\cdot x^m/x_\rho & \text{if}\quad m-e_\rho\in \sigma_M\\
0 & \text{otherwise} 
\end{array}\right .
\]
It follows that 
\[
(Q_\rho)^\sigma_m\cong
\left\{
\begin{array}{ll}
k\cdot x^m & \text{if}\quad m\in \sigma_M\cap \rho^\perp\\
0 & \text{otherwise}   
\end{array}\right .
\]
The set $I_m:=\{\rho \in \sigma(1) \mid m \in \sigma_M\cap \rho^\perp\}$
counts the non-zero graded components $(Q_\rho)^\sigma_m$. Thus the
map $\beta_m^\sigma$ is represented as follows in coordinates
\[
\UseComputerModernTips
\xymatrix{
E_m^\sigma\ar[r] \ar[d]_\approx & {\bigoplus}_{\rho\in I_m}(Q_\rho)^\sigma_m
\ar[d]_\approx \\
k^n \ar[r]^{B_m} & k^{I_m}
}
\]
by the linear map $B_m$ given as 
$$ \lambda\mapsto (\langle\lambda,n(\rho)\rangle)_{\rho \in I_m},$$
see \ref{tordiv}. Now $B_m$ is surjective if and only if the vectors
$n(\rho)$, considered as vectors in $\Z^d\otimes_\Z k,$ 
are linearly independent for $\rho\in I_m$. It follows that $\beta_m^\sigma$
is surjective for any $m\in M$ if and only if all the elements
$\{n(\rho)\}_{\rho \in \sigma(1)}$ are linearly independent or, equivalently,
if and only if $\sigma$ is simplicial.
\end{proof}
\end{sub}
\vskip3mm

The complement of the union of those $U_\sigma$ where $\sigma$ is simplicial
is called the
non-simplicial locus and denoted $S'(X)$. It is easy to see that then
\[
S'(X)=\underset{\sigma\in \Delta'}{\bigcup} \overline{\orb(\sigma)}
\]
where $\Delta'\subset\Delta$ denotes the set of non-simplicial cones. Because
any $\sigma$ of dimension 2 is simplicial, 
$\codim\, \overline{\orb(\sigma)} \geq 3$ for any nonsimplicial $\sigma$ and so
\[
\codim S'(X)\geq 3.
\]
\vskip3mm

\begin{sub}{\bf Corollary:}\label{fitt}
With the notation above 

\begin{enumerate}
\item [(i)]  $\Supp(\kc)= S'(X)$ and $S'(X)$
  is the subvariety of the Fitting ideal of $\alpha$ of its maximal minors
\item [(ii)] $\widetilde{\Omega}^1_{X, p}$ is maximally Cohen--Macaulay if and 
only if
  $p\notin S'(X)$.
\end{enumerate}

\begin{proof}  (i) follows from lemma \ref{sipl} and the diagrams in
\ref{zariskisheaves}.
Because $X$ and each divisor $D_\rho$ is Cohen--Macaulay, also $\ko(-D_\rho)$
is Cohen--Macaulay. It follows that
\[
\text{ depth } \widetilde{\Omega}^1_{X, p} = \text{ depth } \kc_p+2\leq n-1
\]
for any $p\in S'(X)$ because $\codim S'(X)\geq 3$. This proves (ii).
\end{proof}
\end{sub} 
\vskip5mm

\begin{sub}{\bf The primary decomposition of $\widetilde{\Omega}^1_X$}
\label{decdiff}\rm\\
The primary decomposition of $\widetilde{\Omega}^1_X$ in 
$M\otimes _\Z\ko_X$
is now determined by the sequence of $\beta$. Let $\kf_{\rho_0}$ for any
$\rho_0\in\Delta(1)$ denote the kernel of the composition
$M\otimes_\Z\ko_X\xrightarrow{\beta} \underset{\rho}{\oplus}\ko_{D_\rho}\to
\ko_{D_{\rho_0}}$. Then
\[
\widetilde{\Omega}^1_X=\underset{\rho\in\Delta(1)}{\cap}\kf_\rho.
\]
\end{sub}

\begin{sub}{\bf Lemma:}\label{decom}
This is the primary decomposition of $\widetilde{\Omega}^1_X$ in
$M\otimes_\Z\ko_X$.
\begin{proof}
The cokernel $\kc$ in the above Ishida complex is zero outside the
singular (or even non--simplicial) locus of codimension $\geq 2$ by \cite{Od},
Theorem 3.6. Therefore, the homomorphisms $M\otimes_\Z\ko_X \to \ko_{D_\rho}$
are essentially surjective and their images $\kl_\rho$ have $D_\rho$ as
support. Because $D_\rho$ is a toric variety itself, $\ko_{D_\rho}$ and
$\kl_\rho$ don't have torsion of dimension $<n-1$. This means that $\kf_\rho$
is $D_\rho$--primary.
\end{proof}
\end{sub}
\vskip3mm

\begin{sub}{\bf Remark:}\label{aspdec}\quad\rm
The above primary decomposition can be viewed under the aspect of
fine-graded modules and under the aspect of the isotypical decomposition
of the modules of local sections as treated in \cite{P1},\cite{P4}.  

Firstly, a sheaf $\kf_\rho$ is
the sheafification of the fine-graded module
$F_\rho : =\ker(M\otimes S\to S/S(-e_\rho))$,
and 
$$F=\underset{\rho\in\Delta(1)}{\cap}F_\rho$$
is the fine-graded primary decomposition of $F$, inducing that of 
$\widetilde{F}=\widetilde{\Omega}^1_X.$

Secondly, with the notation as in the above proof, the spaces 
$F^\sigma_m\cong \Gamma(U_\sigma, \widetilde{\Omega}^1_X)_m$
are determined as kernels as follows. Let $V:=M\otimes_\Z k$ and 
$V(\rho):=\rho^\perp\otimes_\Z k\subset V.$ Then
\[
F^\sigma_m=
\left\{
\begin{array}{ll}
0 & \text{if}\quad m\not\in \sigma_M\\
\underset{\rho \in I_m}{\cap} V(\rho) & \text{if}\quad m\in \sigma_M
\quad \text{and}\quad I_m\not=\emptyset\\ 
V &\text{if}\quad m\in \sigma_M\quad \text{and}\quad I_m=\emptyset
\end{array}\right .
.\]
Similarly,
\[
(F_\rho)^\sigma_m=
\left\{
\begin{array}{ll}
0 & \text{if}\quad m\not\in \sigma_M\\
V(\rho) & \text{if}\quad m\in \sigma_M\cap \rho^\perp\\
 V &\text{otherwise}
\end{array}\right .
,\]
such that 
$$F^\sigma_m=\underset{\rho\in\Delta(1)}{\cap}(F_\rho)^\sigma_m\ .$$
\end{sub}

We conclude with two examples which demonstrate that a homogeneous
$S$-module can be more complicated than its sheafification if the
quotient $H$ does not act freely on $\hat{X}$.
Hence, graded primary decompositions must not necessarily descent.

\begin{sub}\label{singexample1}\rm {\bf Example:}\;
Consider the quotient of $k^2$ by $\Z / 2\Z$, where $\Z / 2\Z$ acts by
mapping $v$ to $-v$ for every $v \in k^2$. The induced grading on the
polynomial ring $k[x,y]$ is given by $\deg_{\Z / 2\Z}(x) = \deg_{\Z / 2\Z}(y) = 1$.
The sheafification of any $k[x,y]$-module which is concentrated in
degree $1$ to $\Spec(k[x,y]^{\Z / 2\Z})$ is zero.
\end{sub}

\begin{sub}\label{singexample2}\rm {\bf Example:}\;  
Consider the polynomial ring in four variables $S = k[x_1, \dots, x_4]$, which 
is
$\mathbb{Z}$-graded by setting $\deg_A x_1 = \deg_A x_4 = 1$ and $\deg_A x_2 =
\deg_A x_3 = -1$ and $B = S$. Let $I := \langle x_1^2, x_2, x_3, x_4
\rangle$ and let $E$  be the image of the ideal $\langle x_1 \rangle$ in the 
quotient
$S / I$. This is a torsion module whose dimension as $k$-vector space is 1 and
whose only degree is 1. In fact, it is a primary module whose associated prime
ideal is the maximal ideal $\langle x_1, \dots, x_4 \rangle$ of $S$. Taking
degree zero then leads to the zero module.
\end{sub}

\end{document}